\documentclass[12pt]{article}
\usepackage{amsmath,amsfonts,amssymb,amsthm}
\input amssym.def
\begin{document}
\def \Z{\Bbb Z}
\def \C{\Bbb C}
\def \R{\Bbb R}
\def \Q{\Bbb Q}
\def \N{\Bbb N}

\def \A{{\mathcal{A}}}
\def \D{{\mathcal{D}}}
\def \E{{\mathcal{E}}}
\def \H{\mathcal{H}}
\def \S{{\mathcal{S}}}
\def \wt{{\rm wt}}
\def \tr{{\rm tr}}
\def \span{{\rm span}}
\def \Res{{\rm Res}}
\def \Der{{\rm Der}}
\def \End{{\rm End}}
\def \Ind {{\rm Ind}}
\def \Irr {{\rm Irr}}
\def \Aut{{\rm Aut}}
\def \GL{{\rm GL}}
\def \Hom{{\rm Hom}}
\def \mod{{\rm mod}}
\def \ann{{\rm Ann}}
\def \ad{{\rm ad}}
\def \rank{{\rm rank}\;}
\def \<{\langle}
\def \>{\rangle}

\def \g{{\frak{g}}}
\def \h{{\hbar}}
\def \k{{\frak{k}}}
\def \sl{{\frak{sl}}}
\def \gl{{\frak{gl}}}

\def \be{\begin{equation}\label}
\def \ee{\end{equation}}
\def \bex{\begin{example}\label}
\def \eex{\end{example}}
\def \bl{\begin{lem}\label}
\def \el{\end{lem}}
\def \bt{\begin{thm}\label}
\def \et{\end{thm}}
\def \bp{\begin{prop}\label}
\def \ep{\end{prop}}
\def \br{\begin{rem}\label}
\def \er{\end{rem}}
\def \bc{\begin{coro}\label}
\def \ec{\end{coro}}
\def \bd{\begin{de}\label}
\def \ed{\end{de}}

\newcommand{\m}{\bf m}
\newcommand{\n}{\bf n}
\newcommand{\nord}{\mbox{\scriptsize ${\circ\atop\circ}$}}
\newtheorem{thm}{Theorem}[section]
\newtheorem{prop}[thm]{Proposition}
\newtheorem{coro}[thm]{Corollary}
\newtheorem{conj}[thm]{Conjecture}
\newtheorem{example}[thm]{Example}
\newtheorem{lem}[thm]{Lemma}
\newtheorem{rem}[thm]{Remark}
\newtheorem{de}[thm]{Definition}
\newtheorem{hy}[thm]{Hypothesis}
\makeatletter \@addtoreset{equation}{section}
\def\theequation{\thesection.\arabic{equation}}
\makeatother \makeatletter

\begin{center}
{\Large \bf Quantum vertex algebras and their $\phi$-coordinated
quasi modules}
\end{center}

\begin{center}
{Haisheng Li\footnote{Partially supported by NSF grant
DMS-0600189}\\
Department of Mathematical Sciences\\
Rutgers University, Camden, NJ 08102}
\end{center}

\begin{abstract}
We develop a theory of $\phi$-coordinated (quasi) modules for a
nonlocal vertex algebra and we establish a conceptual construction
of nonlocal vertex algebras and their $\phi$-coordinated (quasi)
modules, where $\phi$ is what we call an associate of the
one-dimensional additive formal group. By specializing $\phi$ to a
particular associate, we obtain a new construction of weak quantum
vertex algebras in the sense of \cite{li-qva1}. As an application,
we associate weak quantum vertex algebras to quantum affine
algebras, and we also associate quantum vertex algebras and
$\phi$-coordinated modules to a certain quantum
$\beta\gamma$-system.
\end{abstract}

\section{Introduction}
In the general field of vertex algebras, arguably a central problem
(see \cite{fj}, \cite{efk}) is to develop a suitable theory of
quantum vertex algebras so that quantum vertex algebras can be
associated to quantum affine algebras in the same way that vertex
algebras were associated to affine Lie algebras. With solving this
problem as one of the main goals, in the past we have conducted a
series of studies. In \cite{li-qva1}, we formulated and studied a
notion of (weak) quantum vertex algebra, inspired by
Etingof-Kazhdan's notion of quantum vertex operator algebra (see
\cite{ek}), and we established a conceptual construction of nonlocal
vertex algebras and weak quantum vertex algebras together with their
modules. Nonlocal vertex algebras (which are weak $G_{1}$-vertex
algebras in the sense of \cite{li-g1} and are essentially field
algebras in the sense of \cite{bk}) are analogs of noncommutative
associative algebras, in contrast with that vertex algebras are
analogs of commutative associative algebras. Furthermore, weak
quantum vertex algebras are nonlocal vertex algebras that satisfy a
certain braided locality (commutativity). As an application, we had
associated nonlocal vertex algebras to quantum affine algebras.
Unfortunately, the associated nonlocal vertex algebras are not weak
quantum vertex algebras. A crucial question is whether this theory
of (weak) quantum vertex algebras is the (or a) right one for
solving the aforementioned problem.

The main goal of this paper is to answer this very question. In this
paper, we develop a theory of what we call $\phi$-coordinated
(quasi) modules for nonlocal vertex algebras (including vertex
algebras and weak quantum vertex algebras) and we establish a
conceptual construction of nonlocal vertex algebras and their
$\phi$-coordinated (quasi) modules.  In this new theory, the
parameter $\phi$, which is an element of $\C((x))[[z]]$ satisfying
certain conditions, is what we call an associate of the
one-dimensional additive formal group $F_{\rm a}(x,y)=x+y$.  When
$\phi(x,z)=x+z$ (the additive formal group itself), this
construction of nonlocal vertex algebras reduces to the construction
of \cite{li-qva1} while the notion of $\phi$-coordinated (quasi)
module reduces to the notion of (quasi) module. Specializing $\phi$
to another particular associate of $F_{\rm a}(x,y)$, we obtain a new
construction of weak quantum vertex algebras, which enables us to
associate weak quantum vertex algebras to quantum affine algebras
through $\phi$-coordinated quasi modules.

We now go into some technical details to describe the contents of
this paper. Let $W$ be a general vector space and set $\E(W)=\Hom
(W,W((x)))$. In \cite{li-qva1}, we studied the vertex algebra-like
structures generated by various types of subsets of $\E(W)$, where
the most general type consists of what were called quasi compatible
subsets.  On $\E(W)$, we considered partial operations
$(a(x),b(x))\mapsto a(x)_{n}b(x)$ for any quasi compatible pair
$(a(x),b(x))$ and for $n\in \Z$. Roughly speaking, they were defined
in terms of the generating function $Y_{\E}(a(x),z)b(x)=\sum_{n\in
\Z}a(x)_{n}b(x)z^{-n-1}$ by
$$``Y_{\E}(a(x),z)b(x)=[a(x_{1})b(x)]|_{x_{1}=x+z},"$$
which essentially uses what physicists call the operator product
expansion. It was proved therein that any quasi compatible subset of
$\E(W)$ generates a nonlocal vertex algebra with $W$ as a quasi
module in a certain sense. This generalizes the corresponding result
of \cite{li-g1}. Furthermore, it was proved that every what we
called (resp. quasi) $\S$-local subset of $\E(W)$ generates a weak
quantum vertex algebra with $W$ as a (resp. quasi) module, which
generalizes the corresponding results of \cite{li-local} and
\cite{li-gamma}.

The essence of this present paper is a family generalization of the
vertex operator operation $Y_{\E}$, parameterized by a formal series
$\phi(x,z)\in \C((x))[[z]]$, satisfying
$$\phi(x,0)=x,\ \ \ \ \phi(\phi(x,y),z)=\phi(x,y+z).$$
We call such a formal series an {\em associate} of the
one-dimensional additive formal group $F_{\rm a}(x,y)=x+y$, where a
formal group to an associate is like a group $G$ to a $G$-set. It is
proved that for any $p(x)\in \C((x))$, $e^{zp(x)(d/dx)}x$ is an
associate and that every associate is of this form. In particular,
we have $\phi(x,z)=x+z=F_{\rm a}(x,z)$ for $p(x)=1$ and
$\phi(x,z)=xe^{z}$ for $p(x)=x$. For a quasi compatible pair
$(a(x),b(x))$ in $\E(W)$, we define $a(x)_{n}^{\phi}b(x)$ for $n\in
\Z$ in terms of the generating function
$$Y_{\E}^{\phi}(a(x),z)b(x)=\sum_{n\in \Z}a(x)_{n}^{\phi}b(x)z^{-n-1}$$
roughly by
$$``Y_{\E}^{\phi}(a(x),z)b(x)=[a(x_{1})b(x)]|_{x_{1}=\phi(x,z)}"$$
(see Section 2 for the precise definition). We prove that any quasi
compatible subset $U$ of $\E(W)$ generates a nonlocal vertex algebra
under the operation $Y_{\E}^{\phi}$. To describe the relationship
between such nonlocal vertex algebras and the space $W$, we
introduce a notion of $\phi$-coordinated (quasi) module. In terms of
this notion, the space $W$ becomes a $\phi$-coordinated (quasi)
module for those nonlocal vertex algebras.

To deal with quantum affine algebras, we formulate notions of quasi
$\S_{trig}$-local subset and $\S_{trig}$-local subset of $\E(W)$
with $W$ a general vector space. We prove that every quasi
$\S_{trig}$-local subset generates under $Y_{\E}^{\phi}$ with
$\phi(x,z)=xe^{z}$ a weak quantum vertex algebra in the sense of
\cite{li-qva1}. If $W$ is taken to be a highest weight module for a
quantum affine algebra, the  generating functions of the generators
in the Drinfeld realization form a quasi $\S_{trig}$-local subset
and hence they generate a weak quantum vertex algebra with $W$ as a
$\phi$-coordinated quasi module by our conceptual result. In this
way, we obtain a canonical association of quantum affine algebras
with weak quantum vertex algebras. In a sequel, we shall study this
association in a deeper level, to determine the structure of the
associated weak quantum vertex algebras and prove that they are
indeed quantum vertex algebras.

On the other hand, as a toy example we apply this general machinery
to a certain quantum $\beta\gamma$-system. To this system we
associate quantum vertex algebras and $\phi$-coordinated modules
explicitly. This particular quantum $\beta\gamma$-system is in fact
a one-dimensional trigonometric type Zamolodchikov-Faddeev algebra
(see \cite{fad}, \cite{zam}). Previously, rational type
Zamolodchikov-Faddeev algebras have been associated with quantum
vertex algebras and modules (see \cite{li-qva1}, \cite{li-qva2},
\cite{kl}). The quantum vertex algebras associated to the
trigonometric type quantum $\beta\gamma$-system in this paper are
described by a certain rational type quantum $\beta\gamma$-system.

While this paper has set up a basic foundation for the theory of
$\phi$-coordinated (quasi) modules for nonlocal vertex algebras,
including vertex (operator) algebras in particular, there are many
aspects to be explored further. We note that though this paper is
among a series of papers in a long program, it is pretty much self
contained and can be read independently.

This paper is organized as follows: In Section 2, we define the
notion of associate of the additive formal group and we construct
and classify all the associates. In Section 3, we introduce a notion
of $\phi$-coordinated (quasi) module for a nonlocal vertex algebra.
In Section 4, we give a conceptual construction of nonlocal vertex
algebras and their $\phi$-coordinated (quasi) modules. In Section 5,
we study $\phi$-coordinated modules for (weak) quantum vertex
algebras with $\phi(x,z)=xe^{z}$. In Section 6, we study two quantum
$\beta\gamma$-systems in terms of quantum vertex algebras and their
$\phi$-coordinated modules.

\section{Associates of the one-dimensional additive formal group}
In this section, we formulate and study a notion of associate for a
one-dimensional formal group. For the one-dimensional additive
formal group we construct and classify its associates.

Throughout this paper, we use the usual symbols $\C$ for the complex
numbers, $\Z$ for the integers, and $\N$ for the nonnegative
integers. For this paper, we shall be working on $\C$ and we use the
fairly standard formal variable notations and conventions (see
\cite{flm}, \cite{fhl}; cf. \cite{ll}).

We first recall the notion of formal group (see \cite{boch}).

\bd{dformal-group} {\em  A {\em one-dimensional formal group} over
$\C$ is a formal power series $F(x,y)\in \C[[x,y]]$ such that
\begin{eqnarray*}
&&F(x,y)= x+y +\mbox{terms of high degree},\\
&&F(x,F(y,z))=F(F(x,y),z).
\end{eqnarray*}}
\ed

The simplest example is the one-dimensional {\em additive formal
group}
\begin{eqnarray}
F_{\rm a}(x,y)=x+y.
\end{eqnarray}

We formulate the following notion, which is an analog of the notion
of $G$-set for a group $G$ to a certain extent:

\bd{dassociate} {\em  Let $F(x,y)$ be a one-dimensional formal group
over $\C$. An {\em associate} of $F(x,y)$ is a formal series
$\phi(x,z)\in \C((x))[[z]]$, satisfying the condition that
\begin{eqnarray}\label{eformal-asso}
&&\hspace{1.5cm}\phi(x,0)=x,\nonumber\\
&& \phi(\phi(x,x_{2}),x_{0})=\phi(x,F(x_{0},x_{2})).
\end{eqnarray}}
\ed

\br{rclarification} {\em We here verify the well-definedness of the
two expressions on both sides of (\ref{eformal-asso}). Note that as
$\phi(x,z)\in \C((x))[[z]]$ with $\phi(x,0)=x$, $\phi(x,z)$ is a
unit in the algebra $\C((x))[[z]]$, so that it is well understood
that
\begin{eqnarray}
\phi(x,z)^{m}\in \C((x))[[z]]\ \ \ \mbox{ for }m\in \Z.
\end{eqnarray}
Write $\phi(x,z)=x+zA$ with $A\in \C((x))[[z]]$. For
$f(x)=\sum_{m\ge k}a_{m}x^{m}\in \C((x))$ with $k\in \Z,\; a_{m}\in
\C$, we have
$$f(\phi(x,x_{2}))=\sum_{m\ge k}a_{m}\phi(x,x_{2})^{m}
=\sum_{m\ge k}\sum_{i\ge 0}\binom{m}{i}a_{m}x^{m-i}x_{2}^{i}A^{i}\in
\C((x))[[x_{2}]].$$ Set $\phi(x,z)=\sum_{n\ge 0}f_{n}(x)z^{n}$ with
$f_{n}(x)\in \C((x))$. By definition we have
$$\phi(\phi(x,x_{2}),x_{0})=\sum_{n\ge 0}f_{n}(\phi(x,x_{2}))x_{0}^{n}
\in \C((x))[[x_{0},x_{2}]].$$ On the other hand, as $F(0,0)=0$,
$\phi(x,F(x_{0},x_{2}))$ also exists in  $\C((x))[[x_{0},x_{2}]]$.}
\er

For this paper, our interest is on the additive formal group $F_{\rm
a}(x,y)$. The following is an explicit construction of associates
for $F_{\rm a}(x,y)$:

\bp{pall-assoc} Let $p(x)\in \C((x))$. Set
$$\phi(x,z)=e^{zp(x)\frac{d}{dx}}x
=\sum_{n\ge 0}\frac{z^{n}}{n!}\left(p(x)\frac{d}{dx}\right)^{n}x\in
\C((x))[[z]].$$ Then $\phi(x,z)$ is an associate of $F_{\rm
a}(x,y)$. Furthermore, every associate $\phi(x,z)$ of $F_{\rm
a}(x,y)$ is of this form with $p(x)$ uniquely determined. \ep

\begin{proof} For the first assertion, clearly, $\phi(x,0)=x$.
Since $e^{x_{2}p(x)(d/dx)}$ is an automorphism of the algebra
$\C((x))[[x_{0},x_{2}]]$, we have
$$e^{x_{2}p(x)(d/dx)}f(x,x_{0},x_{2})=f(e^{x_{2}p(x)(d/dx)}x,x_{0},x_{2})$$
for $f(x,x_{0},x_{2})\in \C((x))[[x_{0},x_{2}]]$. Then
\begin{eqnarray*}
&&\phi(x,x_{0}+x_{2})=e^{(x_{0}+x_{2})p(x)(d/dx)}x
=e^{x_{2}p(x)(d/dx)}e^{x_{0}p(x)(d/dx)}x\\
&&\ \ \ \ =e^{x_{2}p(x)(d/dx)}\phi(x,x_{0})
=\phi(e^{x_{2}p(x)(d/dx)}x,x_{0})=\phi(\phi(x,x_{2}),x_{0}).
\end{eqnarray*}
This proves that $\phi(x,z)$ is an associate of $F_{\rm a}(x,y)$.

On the other hand, let $\phi(x,z)$ be any associate. Denote the
formal partial derivatives of $\phi(x,z)$ by $\phi_{x}(x,z)$ and
$\phi_{z}(x,z)$. Set $p(x)=\phi_{z}(x,0)\in \C((x))$. We claim that
$\phi(x,z)=e^{zp(x)(d/dx)}x$. By definition we have
$$\phi(\phi(x,y),z)=\phi(x,y+z).$$
Extracting the coefficients of $z$, we get
$$\phi_{z}(\phi(x,y),0)=\phi_{y}(x,y).$$
Noticing that $\phi_{z}(\phi(x,y),0)=p(\phi(x,y))$, we see that
$\phi(x,z)$ satisfies the differential equation
$$\frac{\partial}{\partial z}\phi(x,z)=p(\phi(x,z))$$
with initial condition $\phi(x,0)=x$. For any $k\in \N$, we have
$$\left(\frac{\partial}{\partial z}\right)^{k+1}\phi(x,z)
=\left(\frac{\partial}{\partial z}\right)^{k}p(\phi(x,z)),$$ which
gives a recursion of the coefficients of $z^{n}$ in $\phi(x,z)$ for
$n\in \N$. It follows that the initial value problem has a unique
solution. Thus $\phi(x,z)=e^{zp(x)(d/dx)}x$. We also see that $p(x)$
is uniquely determined by $p(x)=\phi_{z}(x,0)$.
\end{proof}

\bex{associates} {\em  Here we work out some special examples by
using Proposition \ref{pall-assoc}. We have $\phi(x,z)=x$ for
$p(x)=0$, $\phi(x,z)=e^{z(d/dx)}x=x+z=F_{\rm a}(x,z)$ for $p(x)=1$,
$\phi(x,z)=e^{zx(d/dx)}x=xe^{z}$ for $p(x)=x$. We also have
$\phi(x,z)=e^{zx^{2}(d/dx)}x=\frac{x}{1-zx}$ for $p(x)=x^{2}$ and
$\phi(x,z)=x+\log(1+z/x)$ for $p(x)=x^{-1}$. } \eex

\br{rsubstitution} {\em We here discuss a certain formal
substitution slightly different from those discussed in Remark
\ref{rclarification}. Let $\phi(x,z)\in \C((x))[[z]]$ with
$\phi(x,0)=x$ and let $f(x_{1},x)=\sum_{m,n\ge
k}a(m,n)x_{1}^{m}x^{n}\in \C((x_{1},x))$ with $k\in \Z$. Set
\begin{eqnarray*}
f(\phi(x,z),x)=\sum_{m,n\ge k}a(m,n)\phi(x,z)^{m}x^{n} =\sum_{m,n\ge
k}\sum_{i\ge 0}\binom{m}{i}a(m,n)z^{i}A^{i}x^{m+n-i},
\end{eqnarray*}
which exists in $\C((x))[[z]]$, where $\phi(x,z)=x+zA$ with $A\in
\C((x))[[z]]$. We have
\begin{eqnarray}
f(\phi(x,z),x)
=\Res_{x_{1}}x_{1}^{-1}\delta\left(\frac{\phi(x,z)}{x_{1}}\right)f(x_{1},x).
\end{eqnarray}
Furthermore, for $A(x_{1},x_{2})\in \Hom (W,W((x_{1},x_{2})))$ with
$W$ a vector space, we have
\begin{eqnarray}
A(\phi(x_{2},z),x_{2})\in \Hom (W,W((x_{2}))[[z]])= (\Hom
(W,W((x_{2}))))[[z]].
\end{eqnarray}}
\er

The following technical result plays an important role later:

\bl{lpower-k} Let $\phi(x,z)$ be an associate of $F_{\rm a}(x,y)$
with $\phi(x,z)\ne x$. Then
$$f(\phi(x,z),x)\ne 0\ \ \ \mbox{ for any nonzero }f(x_{1},x)\in \C((x_{1},x)).$$
\el

\begin{proof} It suffices to prove the assertion for $f(x_{1},x)\in \C[[x_{1},x]]$.
Now assume that $f(\phi(x,z),x)=0$ for some
$$f(x_{1},x)=\sum_{m,n\ge 0}a(m,n)x_{1}^{m}x^{n}\in \C[[x_{1},x]].$$
We claim
\begin{eqnarray}\label{eclaim}
\sum_{m,n\ge 0}m^{r}a(m,n)x^{m+n}=0\ \ \ \mbox{ for all }r\ge 0,
\end{eqnarray}
where $m^{r}=1$ for $m=r=0$ as a local convention. By Proposition
\ref{pall-assoc}, we have $\phi(x,z)=e^{zp(x)(d/dx)}x$ for some
nonzero $p(x)\in \C((x))$. As $f(\phi(x,z),x)=0$, we have
\begin{eqnarray}\label{ewehave}
\sum_{m,n\ge 0}a(m,n)x^{n}e^{zp(x)(d/dx)}x^{m} =\sum_{m,n\ge
0}a(m,n)x^{n}\left(e^{zp(x)(d/dx)}x\right)^{m}=0.
\end{eqnarray}
Extracting the constant term with respect to $z$ (equivalently
setting $z=0$), we get
$$\sum_{m,n\ge 0}a(m,n)x^{m+n}=0,$$
proving the base case with $r=0$. Extracting the coefficient of $z$
we get
$$\left(\sum_{m,n\ge 0}ma(m,n)x^{m+n}\right)p(x)x^{-1}=0,$$
which implies
$$\sum_{m,n\ge 0}ma(m,n)x^{m+n}=0,$$
confirming the case with $r=1$. Assume it is true for $0\le r\le k$
with $k\ge 1$. Then
\begin{eqnarray}\label{einductivestep}
\sum_{m,n\ge 0}a(m,n)x^{n}g(x)\left(x\frac{d}{dx}\right)^{r}x^{m}=
g(x)\sum_{m,n\ge 0}m^{r}a(m,n)x^{m+n}=0
\end{eqnarray}
for $0\le
r\le k,\; g(x)\in \C((x))$. Extracting the coefficient of $z^{k+1}$
{}from (\ref{ewehave}) we have
\begin{eqnarray}\label{ek+1}
\sum_{m,n\ge
0}a(m,n)x^{n}\left(p(x)\frac{d}{dx}\right)^{k+1}x^{m}=0.
\end{eqnarray}
 Noticing that
$$\left(p(x)\frac{d}{dx}\right)^{k+1}-(x^{-1}p(x))^{k+1}\left(x\frac{d}{dx}\right)^{k+1}
\in \sum_{j=1}^{k}\C((x))\left(x\frac{d}{dx}\right)^{j},$$ using
(\ref{ek+1}) and (\ref{einductivestep}) we get
$$\sum_{m,n\ge 0}a(m,n)x^{n}(x^{-1}p(x))^{k+1}\left(x\frac{d}{dx}\right)^{k+1}x^{m}=0.$$
Then we obtain
$$\sum_{m,n\ge 0}m^{k+1}a(m,n)x^{m+n}=0,$$
completing the induction and proving (\ref{eclaim}).

Now, {}from (\ref{eclaim}) we get
$$\sum_{m=0}^{l}m^{r}a(m,l-m)=0$$
for all $\ell,r\ge 0$. It follows  that $a(m,l-m)=0$ for all $0\le
m\le l$. Therefore, we have $f(x_{1},x)=0$, proving the assertion.
\end{proof}

\br{rsimplefact} {\em Here, we collect some simple facts which we
need later. Let $\phi(x,z)$ be an associate of $F_{\rm a}(x,y)$. Let
$h(x_{1},x_{0})\in \C((x_{1}))((x_{0}))$. We see that
$$h(\phi(x_{2},x_{0}),x_{0})\ \mbox{ exists in
}\C((x_{2}))((x_{0})).$$ Furthermore, we have
\begin{eqnarray}
\left(h(x_{1},x_{0})|_{x_{1}=\phi(x_{2},x_{0})}\right)|_{x_{2}=\phi(x_{1},-x_{0})}
=h(x_{1},x_{0}),
\end{eqnarray}
noticing that $\phi(\phi(x,-z),z)=\phi(x,0)=x$. For
$A(x_{1},x_{2})\in \C((x_{1},x_{2}))$, we have
\begin{eqnarray}
A(x_{1},x_{2})|_{x_{1}=\phi(x_{2},x_{0})}
=\left(A(x_{1},x_{2})|_{x_{2}=\phi(x_{1},-x_{0})}\right)|_{x_{1}=\phi(x_{2},x_{0})}.
\end{eqnarray}
Assume $\phi(x,z)\ne x$. By Lemma \ref{lpower-k}, for
$A(x_{1},x_{2}),B(x_{1},x_{2})\in \C((x_{1},x_{2}))$, the relation
$$A(x_{1},x_{2})|_{x_{1}=\phi(x_{2},x_{0})}=B(x_{1},x_{2})|_{x_{1}=\phi(x_{2},x_{0})}$$
implies $A(x_{1},x_{2})=B(x_{1},x_{2})$.} \er

\section{$\phi$-coordinated quasi modules for nonlocal vertex algebras}
This is a short preliminary section. In this section, we recall the
definitions of a nonlocal vertex algebra and a (quasi) module, and
we define the notion of $\phi$-coordinated (quasi) module for a
nonlocal vertex algebra. We also give a construction of
$\phi$-coordinated modules through Borcherds's construction of
nonlocal vertex algebras.

We begin by recalling the notion of nonlocal vertex algebra, which
plays a central role in this paper. A {\em nonlocal vertex algebra}
is a vector space $V$ equipped with a linear map
\begin{eqnarray*}
Y(\cdot,x):\ \ V&\rightarrow& \Hom (V,V((x)))\subset (\End V)[[x,x^{-1}]]\\
v&\mapsto& Y(v,x)=\sum_{n\in \Z}v_{n}x^{-n-1}\ \ (\mbox{where
}v_{n}\in \End V),
\end{eqnarray*}
and equipped with a distinguished vector ${\bf 1}\in V$, satisfying
the conditions that
$$Y({\bf 1},x)v=v,\ \ Y(v,x){\bf 1}\in V[[x]]\  \mbox{ and } \ \lim_{x\rightarrow
0}Y(v,x){\bf 1}\ (=v_{-1}{\bf 1})=v\ \ \mbox{ for }v\in V,$$
 and that for $u,v,w\in V$, there exists $l\in
\N$ such that
\begin{eqnarray}
(x_{0}+x_{2})^{l}Y(u,x_{0}+x_{2})Y(v,x_{2})w
=(x_{0}+x_{2})^{l}Y(Y(u,x_{0})v,x_{2})w
\end{eqnarray}
(the {\em weak associativity}).

Let $V$ be a nonlocal vertex algebra, which is fixed throughout this
section. Let $\D$ be the linear operator on $V$ defined by $\D
v=v_{-2}{\bf 1}$ for $v\in V$. Then (see \cite{li-g1})
\begin{eqnarray}
Y(v,x){\bf 1}=e^{x\D}v,\ \ \ \ [\D, Y(v,x)]=Y(\D
v,x)=\frac{d}{dx}Y(v,x).
\end{eqnarray}

Among general nonlocal vertex algebras, what we called weak quantum
vertex algebras in \cite{li-qva1} form a distinguished family. A
{\em weak quantum vertex algebra} is defined by using the same set
of axioms except replacing the weak associativity axiom with the
condition that for any $u,v\in V$, there exist
$$u^{(i)}, v^{(i)}\in V,\;  f_{i}(x)\in \C((x)) \ \ (i=1,\dots,r)$$ such that
\begin{eqnarray}\label{ejacobi-wqva}
&&x_{0}^{-1}\delta\left(\frac{x_{1}-x_{2}}{x_{0}}\right)Y(u,x_{1})Y(v,x_{2})\nonumber\\
&&\hspace{1.5cm}-x_{0}^{-1}\delta\left(\frac{x_{2}-x_{1}}{-x_{0}}\right)
\sum_{i=1}^{r}f_{i}(-x_{0})Y(v^{(i)},x_{2})Y(u^{(i)},x_{1})\nonumber\\
&=&x_{1}^{-1}\delta\left(\frac{x_{2}+x_{0}}{x_{1}}\right)Y(Y(u,x_{0})v,x_{2}).
\end{eqnarray}
A weak quantum vertex algebra can also be defined to be a nonlocal
vertex algebra $V$ that satisfies {\em $\S$-locality}: For any
$u,v\in V$, there exist
$$u^{(i)}, v^{(i)}\in V,\;  f_{i}(x)\in \C((x)) \ \ (i=1,\dots,r)$$ such that
\begin{eqnarray}\label{eslocality-wqva}
(x_{1}-x_{2})^{k}Y(u,x_{1})Y(v,x_{2})
=(x_{1}-x_{2})^{k}\sum_{i=1}^{r}f_{i}(x_{2}-x_{1})Y(v^{(i)},x_{2})Y(u^{(i)},x_{1})
\end{eqnarray}
for some $k\in \N$.

The following notion was introduced in \cite{li-g1}:

\bd{dmodule} {\em A {\em $V$-module} is a vector space $W$ equipped
with a linear map
\begin{eqnarray*}
Y_{W}(\cdot,x):\ \ V&\rightarrow& \Hom (W,W((x)))\subset (\End W)[[x,x^{-1}]]\\
v&\mapsto& Y_{W}(v,x),
\end{eqnarray*}
satisfying the conditions that $Y_{W}({\bf 1},x)=1_{W}$ (the
identity operator on $W$) and that for $u,v\in V,\;w\in W$, there
exists $l\in \N$ such that
\begin{eqnarray}\label{emodule-nva}
(x_{0}+x_{2})^{l}Y_{W}(u,x_{0}+x_{2})Y_{W}(v,x_{2})w
=(x_{0}+x_{2})^{l}Y_{W}(Y(u,x_{0})v,x_{2})w.
\end{eqnarray}} \ed

\br{rsomeproperty}  {\em We note that from \cite{ltw} (Lemma 2.9),
the weak associativity axiom in the definition of a $V$-module can
be equivalently replaced by the condition that for any $u,v\in V$,
there exists $k\in \N$ such that
\begin{eqnarray}\label{econvergence}
(x_{1}-x_{2})^{k}Y_{W}(u,x_{1})Y_{W}(v,x_{2})\in
\Hom(W,W((x_{1},x_{2})))
\end{eqnarray}
and
\begin{eqnarray}\label{eweak-assoc-module}
\left(
(x_{1}-x_{2})^{k}Y_{W}(u,x_{1})Y_{W}(v,x_{2})\right)|_{x_{1}=x_{2}+x_{0}}
=x_{0}^{k}Y_{W}(Y(u,x_{0})v,x_{2}).
\end{eqnarray}
Here, by $\left(
(x_{1}-x_{2})^{k}Y_{W}(u,x_{1})Y_{W}(v,x_{2})\right)|_{x_{1}=x_{2}+x_{0}}$
we mean
$$\iota_{x_{2},x_{0}}\left[\left(
(x_{1}-x_{2})^{k}Y_{W}(u,x_{1})Y_{W}(v,x_{2})\right)|_{x_{1}=x_{2}+x_{0}}\right],$$
the expansion in the nonnegative powers of the second variable
$x_{0}$. Note that
$\left(Y_{W}(u,x_{1})Y_{W}(v,x_{2})\right)|_{x_{1}=x_{2}+x_{0}}$
does {\em not} exist in general. (On the other hand, the
substitution
$\left(Y_{W}(u,x_{1})Y_{W}(v,x_{2})\right)|_{x_{1}=x_{0}+x_{2}}$
{\em always} exists.) Thus  (\ref{econvergence}) is a precondition
for (\ref{eweak-assoc-module}) to make sense. The same principle
also applies to Definitions \ref{dquasimodule} and
\ref{dtrig-module} below.} \er

The following is a modification of the same named notion defined in
\cite{li-qva1}:

\bd{dquasimodule} {\em A {\em quasi $V$-module} is defined as in
Definition \ref{dmodule} except replacing the weak associativity
axiom with the condition that for $u,v\in V$, there exists a nonzero
power series $p(x,y)\in \C[[x,y]]$ such that
$$p(x_{1},x_{2})Y_{W}(u,x_{1})Y_{W}(v,x_{2})\in \Hom
(W,W((x_{1},x_{2}))),$$
\begin{eqnarray}\label{equasi-module-nva}
p(x_{0}+x_{2},x_{2})Y_{W}(Y(u,x_{0})v,x_{2})
=\left(p(x_{1},x_{2})Y_{W}(u,x_{1})Y_{W}(v,x_{2})\right)|_{x_{1}=x_{2}+x_{0}}.
\end{eqnarray}}
\ed

Now, let $\phi(x,z)$ be an associate of the additive formal group
$F_{\rm a}(x,y)$.

\bd{dtrig-module} {\em A {\em $\phi$-coordinated quasi $V$-module}
is defined as in Definition \ref{dmodule} except replacing the weak
associativity axiom with the condition that for $u,v\in V$, there
exists a (nonzero) power series $p(x,y)\in \C[[x,y]]$ such that
$p(\phi(x,z),x)\ne 0$,
\begin{eqnarray}
p(x_{1},x_{2})Y_{W}(u,x_{1})Y_{W}(v,x_{2})\in \Hom
(W,W((x_{1},x_{2}))),
\end{eqnarray}
\begin{eqnarray}
p(\phi(x_{2},x_{0}),x_{2})Y_{W}(Y(u,x_{0})v,x_{2})
=\left(p(x_{1},x_{2})Y_{W}(u,x_{1})Y_{W}(v,x_{2})\right)|_{x_{1}=\phi(x_{2},x_{0})}.
\end{eqnarray}
A {\em $\phi$-coordinated $V$-module} is defined as above except
that $p(x_{1},x_{2})$ is assumed to be a polynomial of the form
$(x_{1}-x_{2})^{k}$ with $k\in \N$.}\ed

\bex{example-class} {\em Let $W$ be a vector space. The space $(\End
W)((x))$ is naturally an associative algebra with identity. Let
$p(x)\in \C((x))$. Then $p(x)\frac{d}{dx}$ is a derivation of $(\End
W)((x))$. By Borcherds' construction, $(\End W)((x))$ becomes a
nonlocal vertex algebra with
$$Y(a(x),z)b(x)=\left(e^{zp(x)\frac{d}{dx}}a(x)\right)b(x)
=a\left(e^{zp(x)\frac{d}{dx}}x\right)b(x)$$ for $a(x),b(x)\in (\End
W)((x))$. Define $Y_{W}(a(x),z)=a(z)$ for $a(x)\in (\End W)((x))$.
We have $Y_{W}(1,z)=1_{W}$, and for $a(x),b(x)\in (\End W)((x))$,
$$Y_{W}(a(x),x_{1})Y_{W}(b(x),x_{2})=a(x_{1})b(x_{2})
\in \Hom (W,W((x_{1},x_{2})))$$ and
\begin{eqnarray*}
& &Y_{W}(Y(a(x),x_{0})b(x),x_{2})=Y(a(x_{2}),x_{0})b(x_{2})
=\left(a(x_{1})b(x_{2})\right)|_{x_{1}=e^{x_{0}p(x_{2})(d/dx_{2})}x_{2}}\\
&&\hspace{2cm}
=\left(Y_{W}(a(x),x_{1})Y_{W}(b(x),x_{2})\right)|_{x_{1}=e^{x_{0}p(x_{2})(d/dx_{2})}x_{2}}.
\end{eqnarray*}
Thus $W$ is a $\phi$-coordinated module for the nonlocal vertex
algebra $(\End W)((x))$ with $\phi(x,z)=e^{zp(x)(d/dx)}x$ and with
$Y_{W}(a(x),x_{0})=a(x_{0})$ for $a(x)\in (\End W)((x))$. In the
next section we shall give a much more sophisticated construction of
nonlocal vertex algebras and their $\phi$-coordinated (quasi)
modules.} \eex

The following is a convenient technical result:

\bl{lbasic-property} Let $V$ be a nonlocal vertex algebra and let
$(W,Y_{W})$ be a $\phi$-coordinated quasi $V$-module. Let $u,v\in V$
and suppose that $q(x_{1},x_{2})\in \C[[x_{1},x_{2}]]$ satisfies
\begin{eqnarray*}
q(x_{1},x_{2})Y_{W}(u,x_{1})Y_{W}(v,x_{2})\in \Hom
(W,W((x_{1},x_{2}))).
\end{eqnarray*}
Then
\begin{eqnarray*}
q(\phi(x_{2},x_{0}),x_{2})Y_{W}(Y(u,x_{0})v,x_{2})
=\left(q(x_{1},x_{2})Y_{W}(u,x_{1})Y_{W}(v,x_{2})\right)|_{x_{1}=\phi(x_{2},x_{0})}.
\end{eqnarray*}
\el

\begin{proof} By definition there exists $p(x_{1},x_{2})\in \C[[x_{1},x_{2}]]$ with
$p(\phi(x,z),x)\ne 0$, satisfying the condition in Definition
\ref{dtrig-module}. Then
\begin{eqnarray*}
&&p(\phi(x_{2},x_{0}),x_{2})q(\phi(x_{2},x_{0}),x_{2})Y_{W}(Y(u,x_{0})v,x_{2})\\
&=&q(\phi(x_{2},x_{0}),x_{2})
\left(p(x_{1},x_{2})Y_{W}(u,x_{1})Y_{W}(v,x_{2})\right)|_{x_{1}=\phi(x_{2},x_{0})}\\
&=&\left(p(x_{1},x_{2})q(x_{1},x_{2})Y_{W}(u,x_{1})Y_{W}(v,x_{2})\right)|_{x_{1}=\phi(x_{2},x_{0})}\\
&=&p(\phi(x_{2},x_{0}),x_{2})
\left(q(x_{1},x_{2})Y_{W}(u,x_{1})Y_{W}(v,x_{2})\right)|_{x_{1}=\phi(x_{2},x_{0})}.
\end{eqnarray*}
Noticing that the powers of $x_{0}$ in both
$$q(\phi(x_{2},x_{0}),x_{2})Y_{W}(Y(u,x_{0})v,x_{2})\ \ \mbox{ and
}\ \
\left(q(x_{1},x_{2})Y_{W}(u,x_{1})Y_{W}(v,x_{2})\right)|_{x_{1}=\phi(x_{2},x_{0})}$$
are truncated from below and that $p(\phi(x_{2},x_{0}),x_{2})\in
\C((x_{2}))((x_{0}))$ is nonzero, we obtain the desired relation by
cancellation.
\end{proof}

We also have the following result:

\bl{lDproperty} Let $V$ be a nonlocal vertex algebra and let
$(W,Y_{W})$ be a $\phi$-coordinated quasi $V$-module. Then
\begin{eqnarray}\label{edproperty}
Y_{W}(e^{x_{0}\D}v,x)=Y_{W}(v,\phi(x,x_{0})) \ \ \ \mbox{ for }v\in
V.
\end{eqnarray}
\el

\begin{proof} For $v\in V$, as $Y_{W}({\bf 1},x)=1_{W}$, we have
$$Y_{W}(v,x_{1})Y_{W}({\bf 1},x_{2})=Y_{W}(v,x_{1})\in
\Hom(W,W((x_{1},x_{2}))).$$
By Lemma \ref{lbasic-property}, we have
\begin{eqnarray*}
Y_{W}(Y(v,x_{0}){\bf 1},x_{2})
=Y_{W}(v,x_{1})|_{x_{1}=\phi(x_{2},x_{0})}=Y_{W}(v,\phi(x_{2},x_{0})),
\end{eqnarray*}
which gives (\ref{edproperty}) as $Y(v,x_{0}){\bf 1}=e^{x_{0}\D}v$.
\end{proof}

\br{rzero-case} {\em Recall from Example \ref{associates} that
$\phi(x,z)=x$ is a particular associate. Let $V$ be a nonlocal
vertex algebra and suppose that $V$ admits a faithful
$\phi$-coordinated quasi module $(W,Y_{W})$ with $\phi(x,z)=x$.  Let
$u,v\in V$. By definition, there exists $p(x_{1},x_{2})\in
\C[[x_{1},x_{2}]]$ with $p(x_{2},x_{2})\ne 0$ such that
$$p(x_{1},x_{2})Y_{W}(u,x_{1})Y_{W}(v,x_{2})\in \Hom
(W,W((x_{1},x_{2}))),$$
$$p(x_{2},x_{2})Y_{W}(Y(u,x_{0})v,x_{2})
=\left(p(x_{1},x_{2})Y_{W}(u,x_{1})Y_{W}(v,x_{2})\right)|_{x_{1}=x_{2}}.$$
Notice that the right-hand side is independent of $x_{0}$. As the
map $Y_{W}$ is assumed to be injective, it follows that $u_{n}v=0$
whenever $n\ne -1$. Consequently, $V$ is merely an ordinary
associative algebra. Now, let $(W,Y_{W})$ be a general
$\phi$-coordinated quasi $V$-module with $\phi(x,z)=x$. It can be
readily seen that $\ker Y_{W}$ is a two-sided ideal of $V$. We have
that $V/(\ker Y_{W})$ is an ordinary associative algebra.} \er

\section{A construction of nonlocal vertex algebras and
their $\phi$-coordinated modules} In this section we present the
conceptual construction of nonlocal vertex algebras and their
$\phi$-coordinated (quasi) modules, by using (quasi) compatible
subsets of formal vertex operators on a vector space. This
generalizes significantly the construction of nonlocal vertex
algebras and their (quasi) modules in \cite{li-g1} and
\cite{li-qva1}.

We begin with certain generalized iota maps (cf. \cite{fhl}).
Denote by $F(\C[[x_{1},x_{2}]])$ the fraction field of the ring
$\C[[x_{1},x_{2}]]$. Since $\C[[x_{1},x_{2}]]$ is also a subring of
the field $\C((x_{1}))((x_{2}))$, there exists a field embedding
\begin{eqnarray}
\iota_{x_{1},x_{2}}: F(\C[[x_{1},x_{2}]])\rightarrow
\C((x_{1}))((x_{2})),
\end{eqnarray}
which is uniquely determined by the condition
$\iota_{x_{1},x_{2}}|_{\C[[x_{1},x_{2}]]}=1$. In fact, we have
$$\iota_{x_{1},x_{2}}|_{\C((x_{1},x_{2}))}=1,$$
noticing that $\C((x_{1},x_{2}))$ is a subalgebra of both $
F(\C[[x_{1},x_{2}]])$ and $\C((x_{1}))((x_{2}))$. This map
$\iota_{x_{1},x_{2}}$ naturally extends the algebra embedding
$\iota_{x_{1},x_{2}}$ of $\C_{*}(x_{1},x_{2})$ into
$\C((x_{1}))((x_{2}))$, which was used in \cite{li-qva1} (cf.
\cite{li-gamma}), where $\C_{*}(x_{1},x_{2})$ is the algebra
extension of $\C[[x_{1},x_{2}]]$ by inverting all the nonzero
polynomials.

\br{rcancelation} {\em We here discuss certain cancellation rules
which shall use extensively in this work. Let $W$ be a vector space.
The space $\Hom (W,W((x)))$ is naturally a vector space over the
field $\C((x))$. Furthermore, $\Hom (W,W((x_{1}))((x_{2})))$ is a
vector space over the field $\C((x_{1}))((x_{2}))$, while $\Hom
(W,W((x_{2}))((x_{1})))$ is a vector space over
$\C((x_{2}))((x_{1}))$. In view of this, for any
$$A(x_{1},x_{2}),B(x_{1},x_{2})\in \Hom (W,W((x_{1}))((x_{2}))),$$
if
$$f(x_{1},x_{2})A(x_{1},x_{2})=f(x_{1},x_{2})B(x_{1},x_{2})$$
for some nonzero $f(x_{1},x_{2})\in \C[[x_{1},x_{2}]]\subset
\C((x_{1}))((x_{2}))$, then $A(x_{1},x_{2})=B(x_{1},x_{2})$. On the
other hand, we have
$$\C((x_{1}))((x_{2}))\cap \C((x_{2}))((x_{1}))=\C((x_{1},x_{2})),$$
and both $\Hom (W,W((x_{1}))((x_{2})))$ and $\Hom
(W,W((x_{2}))((x_{1})))$ are $\C((x_{1},x_{2}))$-modules. (Note that
$\C((x_{1},x_{2}))$ is an algebra but not a field.) In view of this,
for any
$$A(x_{1},x_{2}),B(x_{1},x_{2})\in \Hom (W,W((x_{1}))((x_{2})))+ \Hom
(W,W((x_{2}))((x_{1})))$$ ($\subset (\End W)[[x_{1}^{\pm
1},x_{2}^{\pm 1}]]$), if
$$f(x_{1},x_{2})A(x_{1},x_{2})=f(x_{1},x_{2})B(x_{1},x_{2})$$ for
some invertible element $f(x_{1},x_{2})$ of $\C((x_{1},x_{2}))$,
then $A(x_{1},x_{2})=B(x_{1},x_{2})$.} \er

Let $W$ be a vector space (over $\C$), which is fixed throughout
this section. Set
$$\E(W)=\Hom (W,W((x)))\subset (\End W)[[x,x^{-1}]].$$
The identity operator on $W$ is denoted by $1_{W}$, which is a
typical element of $\E(W)$. Recall the notion of compatibility from
\cite{li-g1}: A finite sequence $a_{1}(x),\dots,a_{r}(x)$ in $\E(W)$
is said to be {\em compatible} if there exists a nonnegative integer
$k$ such that
\begin{eqnarray}
\left(\prod_{1\le i<j\le
r}(x_{i}-x_{j})^{k}\right)a_{1}(x_{1})\cdots a_{r}(x_{r})\in \Hom
(W,W((x_{1},\dots,x_{r}))).
\end{eqnarray}
Furthermore, a subset $U$ of $\E(W)$ is said to be {\em compatible}
if every finite sequence in $U$ is compatible.

We formulate the following notion of quasi compatibility:

\bd{dquasicompatible-generalized} {\em  A finite sequence
$a_{1}(x),\dots,a_{r}(x)$ in $\E(W)$ is said to be {\em quasi
compatible} if there exists a nonzero power series
$p(x_{1},x_{2})\in \C[[x_{1},x_{2}]]$ such that
\begin{eqnarray}\label{epab-condition-phi}
\left(\prod_{1\le i<j\le r}p(x_{i},x_{j})\right)a_{1}(x_{1})\cdots
a_{r}(x_{r})\in \Hom (W,W((x_{1},\dots,x_{r}))).
\end{eqnarray}
Furthermore, a subset $U$ of $\E(W)$ is said to be {\em quasi
compatible} if every finite sequence in $U$ is quasi compatible. }
\ed

Note that this notion of quasi compatibility slightly generalizes
the same named notion defined in \cite{li-qva1} in the way that the
modifier of $p(x_{1},x_{2})$ is changed {}from a nonzero polynomial
to a nonzero power series.

\bd{dphi-quasicompatible} {\em  Let $\phi$ be an associate of
$F_{\rm a}$. A finite sequence $a_{1}(x),\dots,a_{r}(x)$ in $\E(W)$
is said to be {\em $\phi$-quasi compatible} if there exists a power
series $p(x,y)\in \C[[x,y]]$ with $p(\phi(x,z),x)\ne 0$ such that
(\ref{epab-condition-phi}) holds. A subset $U$ of $\E(W)$ is said to
be {\em $\phi$-quasi compatible} if every finite sequence in $U$ is
$\phi$-quasi compatible. } \ed

{}From definition, $\phi$-quasi compatibility implies quasi
compatibility. In fact, for almost all $\phi$, $\phi$-quasi
compatibility is the same as quasi compatibility. We have the
following lemma, where the first assertion is straightforward while
the second immediately follows from Lemma \ref{lpower-k}:

\bl{lphi-quasicompatible}  Let $\phi$ be an associate of $F_{\rm a}$
with $\phi(x,z)\ne x$. Every compatible subset of $\E(W)$ is
$\phi$-quasi compatible and every quasi compatible subset is
$\phi$-quasi compatible. \el

{}From now on, we fix an associate $\phi(x,z)$ of $F_{a}(x,y)$.

\bd{daction} {\em Let $(a(x),b(x))$ be a $\phi$-quasi compatible
pair in $\E(W)$. We define
$$a(x)_{n}^{\phi}b(x)\in \E(W)\ \ \ \mbox{ for }n\in \Z$$
in terms of the generating function
\begin{eqnarray}
Y_{\E}^{\phi}(a(x),z)b(x)=\sum_{n\in \Z}a(x)_{n}^{\phi}b(x) z^{-n-1}
\end{eqnarray}
by
\begin{eqnarray}
Y_{\E}^{\phi}(a(x),z)b(x) =p(\phi(x,z),x)^{-1}
\left(p(x_{1},x)a(x_{1})b(x)\right)|_{x_{1}=\phi(x,z)},
\end{eqnarray}
which lies in $(\Hom (W,W((x))))((z))=\E(W)((z))$, where
$p(x_{1},x_{2})$ is any power series with $p(\phi(x,z),x)\ne 0$ such
that (\ref{epab-condition-phi}) holds and where
$p(\phi(x,z),x)^{-1}$ stands for the inverse of $p(\phi(x,z),x)$ in
$\C((x))((z))$.} \ed

Just as with $Y_{\E}$ (see \cite{li-gamma}), it is straightforward
to show that $Y_{\E}^{\phi}$ is well defined, i.e., the expression
on the right-hand side does not depend on the choice of
$p(x_{1},x_{2})$. {}From definition we have
\begin{eqnarray}
p(\phi(x,z),x)Y_{\E}^{\phi}(a(x),z)b(x)
=\left(p(x_{1},x)a(x_{1})b(x)\right)|_{x_{1}=\phi(x,z)}
\end{eqnarray}
for any power series $p(x_{1},x_{2})$ with $p(\phi(x,z),x)\ne 0$
such that (\ref{epab-condition-phi}) holds.

We shall need the following technical result:

\bl{lproof-need} Let $(a_{i}(x),b_{i}(x))$ $(i=1,\dots,n)$ be
$\phi$-quasi compatible ordered pairs in $\E(W)$. Suppose that
\begin{eqnarray}\label{esum-lemma}
\sum_{i=1}^{n}g_{i}(x_{1},x_{2})a_{i}(x_{1})b_{i}(x_{2})\in \Hom
(W,W((x_{1},x_{2})))
\end{eqnarray}
with $g_{1}(x_{1},x_{2}),\dots,g_{n}(x_{1},x_{2})\in
\C((x_{1},x_{2}))$. Then
\begin{eqnarray}
\sum_{i=1}^{n}g_{i}(\phi(x,z),x)Y_{\E}^{\phi}(a_{i}(x),z)b_{i}(x)
=\left(\sum_{i=1}^{n}g_{i}(x_{1},x)a_{i}(x_{1})b_{i}(x)\right)|_{x_{1}=\phi(x,z)}.
\end{eqnarray}
\el

\begin{proof} There exists $g(x_{1},x_{2})\in \C[[x_{1},x_{2}]]$
with $g(\phi(x,z),x)\ne 0$ such that
$$g(x_{1},x_{2})a_{i}(x_{1})b_{i}(x_{2})\in \Hom (W,W((x_{1},x_{2})))
\;\;\;\mbox{ for }i=1,\dots,n.$$ Then
$$g(\phi(x,z),x)Y_{\E}^{\phi}(a_{i}(x),z)b_{i}(x)
=\left(g(x_{1},x)a_{i}(x_{1})b_{i}(x)\right)\mid_{x_{1}=\phi(x,z)}$$
for $i=1,\dots,n$. Using (\ref{esum-lemma}) we get
\begin{eqnarray*}
&&g(\phi(x,z),x)\sum_{i=1}^{n}g_{i}(\phi(x,z),x)Y_{\E}^{\phi}(a_{i}(x),z)b_{i}(x)
\nonumber\\
&=&\sum_{i=1}^{n}g_{i}(\phi(x,z),x)\left(g(x_{1},x)a_{i}(x_{1})b_{i}(x)\right)|_{x_{1}=\phi(x,z)}
\nonumber\\
&=&\left(g(x_{1},x)\sum_{i=1}^{n}g_{i}(x_{1},x)a_{i}(x_{1})b_{i}(x)\right)|_{x_{1}=\phi(x,z)}
\nonumber\\
&=&g(\phi(x,z),x)\left(\sum_{i=1}^{n}g_{i}(x_{1},x)a_{i}(x_{1})b_{i}(x)\right)|_{x_{1}=\phi(x,z)}.
\end{eqnarray*}
Notice that both $\sum_{i=1}^{n}g_{i}(\phi(x,z),x)
Y_{\E}^{\phi}(a_{i}(x),z)b_{i}(x)$ and
$$\left(\sum_{i=1}^{n}g_{i}(x_{1},x)a_{i}(x_{1})b_{i}(x)\right)|_{x_{1}=\phi(x,z)}$$
lie in $(\Hom (W,W((x)))((z))$. Now, it follows immediately from
cancellation as $g(\phi(x,z),x)\in \C((x))[[z]]$ is nonzero (recall
Remark \ref{rcancelation}).
\end{proof}

\bd{dclosed-Yphi} {\em Let $U$ be a subspace of $\E(W)$ such that
every ordered pair in $U$ is $\phi$-quasi compatible. We say that
$U$ is $Y_{\E}^{\phi}$-{\em closed} if
\begin{eqnarray}
a(x)_{n}^{\phi}b(x)\in U\;\;\;\mbox{ for }a(x),b(x)\in U,\; n\in \Z.
\end{eqnarray}}
\ed

We are going to prove that every $Y_{\E}^{\phi}$-closed $\phi$-quasi
compatible subspace of $\E(W)$, which contains $1_{W}$, is a
nonlocal vertex algebra. First we have:

\bl{lclosed} Assume that $V$ is a subspace of $\E(W)$ such that
every sequence of length $2$ or $3$ in $V$ is $\phi$-quasi
compatible and such that $V$ is $Y_{\E}^{\phi}$-closed. Let
$a(x),b(x),c(x)\in V$ and let $f(x,y)$ be a nonzero power series
such that
\begin{eqnarray}
& &f(y,z)b(y)c(z)\in \Hom (W,W((y,z))),
\label{e4.31}\\
& &f(x,y)f(x,z)f(y,z) a(x)b(y)c(z)\in \Hom
(W,W((x,y,z))).\label{e4.32}
\end{eqnarray}
Then
\begin{eqnarray}
&&f(\phi(x,x_{1}),x)f(\phi(x,x_{2}),x)f(\phi(x,x_{1}),\phi(x,x_{2}))
Y_{\E}^{\phi}(a(x),x_{1})Y_{\E}^{\phi}(b(x),x_{2})c(x)\nonumber\\
&&\ \ \ \ =\left(f(y,x)f(z,x)f(y,z) a(y)b
(z)c(x)\right)|_{y=\phi(x,x_{1}),z=\phi(x,x_{2})}.
\end{eqnarray}
\el

\begin{proof} With (\ref{e4.31}), by Lemma
\ref{lproof-need} we have
\begin{eqnarray}\label{ehphi-theta}
f(\phi(x,x_{2}),x)Y_{\cal{E}}^{\phi}(b(x),x_{2})c(x)
=\left(f(z,x)b(z)c(x)\right)|_{z=\phi(x,x_{2})},
\end{eqnarray}
which gives
\begin{eqnarray}\label{erighthand}
& &f(y,x)f(y,\phi(x,x_{2}))f(\phi(x,x_{2}),x)a(y)
Y_{\cal{E}}^{\phi}(b(x),x_{2})c(x)\nonumber\\
&=&\left(f(y,x)f(y,z) f(z,x) a(y)b(z)c(x)\right)|_{z=\phi(x,x_{2})}.
\end{eqnarray}
{}From (\ref{e4.32}) we see that the right-hand side of
(\ref{erighthand}) lies in $(\Hom (W,W((x,y)))[[x_{2}]]$, so does
the left-hand side. That is,
\begin{eqnarray*}
f(y,x)f(y,\phi(x,x_{2}))f(\phi(x,x_{2}),x)a(y)
Y_{\cal{E}}^{\phi}(b(x),x_{2})c(x) \in (\Hom (W,W((y,x)))[[x_{2}]].
\end{eqnarray*}
Notice that because $b(x)_{m}^{\phi}c(x)=0$ for $m$ sufficiently
large, for every $n\in \Z$, the coefficient of $x_{2}^{n}$ is of the
form
$$\sum_{j=r}^{s}g_{j}(y,x)a(y)(b(x)_{j}^{\phi}c(x)) $$
with $r,s\in \Z$ and $g_{j}(y,x)\in \C((x,y))$. By considering the
coefficient of each power of $x_{2}$ and then using Lemma
\ref{lproof-need}, we have
\begin{eqnarray}
&
&f(\phi(x,x_{1}),x)f(\phi(x,x_{1}),\phi(x,x_{2}))f(\phi(x,x_{2}),x)
Y_{\cal{E}}^{\phi}(a(x),x_{1})Y_{\cal{E}}^{\phi}(b(x),x_{2})c(x)
\nonumber\\
&=&\left(f(y,x)f(y,\phi(x,x_{2}))f(\phi(x,x_{2}),x)
a(y)Y_{\cal{E}}^{\phi}(b(x),x_{2})c(x)\right)|_{y=\phi(x,x_{1})}.
\end{eqnarray}
Using this and (\ref{ehphi-theta}) we obtain
\begin{eqnarray*}
&&f(\phi(x,x_{1}),x)f(\phi(x,x_{1}),\phi(x,x_{2}))f(\phi(x,x_{2}),x)
Y_{\cal{E}}^{\phi}(a(x),x_{1})Y_{\cal{E}}^{\phi}(b(x),x_{2})c(x)\nonumber\\
&=&\left(f(\phi(x,x_{2}),x))f(y,x)f(y,\phi(x,x_{2}))
a(y)Y_{\cal{E}}^{\phi}(b(x),x_{2})c(x)\right)|_{y=\phi(x,x_{1})}\nonumber\\
&=&\left(f(z,x)f(y,x)f(y,z) a(y)b
(z)c(x)\right)|_{y=\phi(x,x_{1}),z=\phi(x,x_{2})},
\end{eqnarray*}
completing the proof.
\end{proof}

Now we are in a position to prove our first key result:

\bt{tclosed} Let $V$ be a subspace of $\E(W)$, that contains
$1_{W}$, such that every sequence of length $2$ or $3$ in $V$ is
(resp. $\phi$-quasi) compatible and $V$ is $Y_{\E}^{\phi}$-closed.
Then $(V,Y_{\E}^{\phi},1_{W})$ carries the structure of a nonlocal
vertex algebra and $W$ is a faithful $\phi$-coordinated (resp.
quasi) $V$-module with $Y_{W}(\alpha(x),x_{0})=\alpha(x_{0})$ for
$\alpha(x)\in V$. \et

\begin{proof} For any $a(x),b(x)\in V$, from definition we have
$a(x)_{n}^{\phi}b(x)=0$ for $n$ sufficiently large and
$a(x)_{n}^{\phi}b(x)\in V$ for any $n\in \Z$ by assumption. We also
have
$$Y_{\E}^{\phi}(1_{W},z)b(x)
=\left(1_{W}(x_{1})b(x)\right)|_{x_{1}=\phi(x,z)}=b(x)$$ and
$$Y_{\E}^{\phi}(a(x),z)1_{W}=\left(a(x_{1})1_{W}\right)|_{x_{1}=\phi(x,z)}=a(\phi(x,z)).$$
Since $\phi(x,z)\in\C((x))[[z]]$ and $\phi(x,0)=x$, we have
$$Y_{\E}^{\phi}(a(x),z)1_{W}\in \E(W)[[z]] \ \
\mbox{ and }\ \ \lim_{z\rightarrow
0}Y_{\E}^{\phi}(a(x),z)1_{W}=a(x).$$
Now, for the assertion on the
nonlocal vertex algebra structure, it remains to prove weak
associativity, i.e., for $a(x),b(x),c(x)\in V$, there exists a
nonnegative integer $k$ such that
\begin{eqnarray*}\label{eweakassocmainthem}
& &(x_{0}+x_{2})^{k}Y_{\cal{E}}^{\phi}(a(x),x_{0}+x_{2})
Y_{\cal{E}}^{\phi}(b(x),x_{2})c(x)\nonumber\\
&=&(x_{0}+x_{2})^{k}Y_{\cal{E}}^{\phi}(Y_{\cal{E}}^{\phi}(a(x),x_{0})b(x),x_{2})c(x).
\end{eqnarray*}
Let $f(x,y)$ be a power series such that $f(\phi(x,z),x)\ne 0$, and
\begin{eqnarray*}
& &f(x,y)a(x)b(y)\in \Hom (W,W((x,y))),\\
& &f(x,y)b(x)c(y)\in \Hom (W,W((x,y))),\\
& &f(x,y)f(x,z)f(y,z) a(x)b(y)c(z)\in \Hom (W,W((x,y,z))).
\end{eqnarray*}
By Lemma \ref{lclosed}, we have
\begin{eqnarray}\label{e5.76}
& &f(\phi(x,x_{2}),x)f(\phi(x,x_{0}+x_{2}),x)f(\phi(x,x_{0}),x)
Y_{\cal{E}}^{\phi}(a(x),x_{0}+x_{2})Y_{\cal{E}}^{\phi}(b(x),x_{2})c(x)
\nonumber\\
&&\ \ \ \ \ =\left(f(z,x)f(y,x)f(y,z)
a(y)b(z)c(x)\right)|_{y=\phi(x,x_{0}+x_{2}),z=\phi(x,x_{2})}.
\end{eqnarray}

On the other hand, let $n\in \Z$ be {\em arbitrarily fixed}. Since
$a(x)_{m}^{\phi}b(x)=0$ for $m$ sufficiently large, there exists a
power series $p(x,y)$, {\em depending on $n$}, such that
$p(\phi(x,z),x)\ne 0$, and
\begin{eqnarray}\label{eges}
p(\phi(x,x_{2}),x) (Y_{\cal{E}}^{\phi}(a(x)_{m}^{\phi}b(x),
x_{2})c(x)
=\left(p(z,x)(a(z)_{m}^{\phi}b(z))c(x)\right)|_{z=\phi(x,x_{2})}
\end{eqnarray}
for {\em all} $m\ge n$. With $f(x,y)a(x)b(y)\in \Hom (W,W((x,y)))$,
we have
\begin{eqnarray}\label{ethis}
f(\phi(x,x_{0}),x)(Y_{\cal{E}}^{\phi}(a(x_{2}),x_{0})b(x_{2}))c(x)
=\left(f(y,x_{2})a(y)b(x_{2})c(x)\right)|_{y=\phi(x_{2},x_{0})}.
\end{eqnarray}
Set
\begin{eqnarray*}
Y_{\cal{E}}^{\phi}(a(x),x_{0})_{\ge n}b(x)=\sum_{m\ge
n}a(x)^{\phi}_{m}b(x).
\end{eqnarray*}
Then for any $q(x)\in \C[[x]]$ we have
\begin{eqnarray}\label{esimplefactresidule}
\Res_{x}x^{n}q(x)Y_{\cal{E}}^{\phi}(a(x),x_{0})b(x)
=\Res_{x}x^{n}q(x)Y_{\cal{E}}^{\phi}(a(x),x_{0})_{\ge n}b(x).
\end{eqnarray}
Using (\ref{esimplefactresidule}), (\ref{eges}) and (\ref{ethis}) we
get
\begin{eqnarray}\label{e5.78}
& &\Res_{x_{0}}x_{0}^{n}f(\phi(x,x_{0}+x_{2}),x)f(\phi(x,x_{0}),x)
p(\phi(x,x_{2}),x) Y_{\cal{E}}^{\phi}(Y_{\cal{E}}^{\phi}(a(x),x_{0})b(x), x_{2})c(x)\nonumber\\
&=&\Res_{x_{0}}x_{0}^{n}f(\phi(x,x_{0}+x_{2}),x)
f(\phi(x,x_{0}),x)p(\phi(x,x_{2}),x)
Y_{\cal{E}}^{\phi}(Y_{\cal{E}}^{\phi}(a(x),x_{0})_{\ge n}b(x), x_{2})c(x)\nonumber\\
&=&\Res_{x_{0}}x_{0}^{n} f(\phi(x,x_{0}+x_{2}),x)f(\phi(x,x_{0}),x)
\left(p(z,x)
Y_{\cal{E}}^{\phi}(a(z),x_{0})_{\ge n}b(z))c(x)\right)|_{z=\phi(x,x_{2})}\nonumber\\
&=&\Res_{x_{0}}x_{0}^{n} f(\phi(x,x_{0}+x_{2}),x)f(\phi(x,x_{0}),x)
 \left(p(z,x)Y_{\cal{E}}^{\phi}(a(z),x_{0})b(z))c(x)\right)|_{z=\phi(x,x_{2})}\nonumber\\
&=&\Res_{x_{0}}x_{0}^{n}\left(f(y,x)f(y,z)
p(z,x)a(y)b(z)c(x)\right)|_{y=\phi(z,x_{0}),z=\phi(x,x_{2})}
\nonumber\\
&=&\Res_{x_{0}}x_{0}^{n}\left(f(y,x)f(y,z)
p(z,x)a(y)b(z)c(x)\right)|_{y=\phi(\phi(x,x_{2}),x_{0}),z=\phi(x,x_{2})}.
\end{eqnarray}
As $\phi(\phi(x,y),z)=\phi(x,y+z)$, combining (\ref{e5.78}) with
(\ref{e5.76}) we get
\begin{eqnarray}\label{e5.79}
& &\Res_{x_{0}}x_{0}^{n}
f(\phi(x,x_{2}),x)f(\phi(x,x_{0}+x_{2}),x)f(\phi(x,x_{0}),x)
 p(\phi(x,x_{2}),x)\nonumber\\
 & &\ \ \ \ \cdot Y_{\cal{E}}^{\phi}(a(x),x_{0}+x_{2})
Y_{\cal{E}}^{\phi}(b(x),x_{2})c(x)\nonumber\\
&=&\Res_{x_{0}}x_{0}^{n}f(\phi(x,x_{2}),x)f(\phi(x,x_{0}+x_{2}),x)f(\phi(x,x_{0}),x)
 p(\phi(x,x_{2}),x)\nonumber\\
& &\ \ \ \ \cdot
Y_{\cal{E}}^{\phi}(Y_{\cal{E}}^{\phi}(a(x),x_{0})b(x), x_{2})c(x).
\end{eqnarray}
Notice that both sides of (\ref{e5.79}) involve only finitely many
negative powers of $x_{2}$. Multiplying both sides by
$p(\phi(x,x_{2}),x)^{-1}f(\phi(x,x_{2}),x)^{-1}$ $(\in
\C((x))((x_{2})))$ we get
\begin{eqnarray*}
& &\Res_{x_{0}}x_{0}^{n}f(\phi(x,x_{0}+x_{2}),x) f(\phi(x,x_{0}),x)
Y_{\cal{E}}^{\phi}(a(x),x_{0}+x_{2})Y_{\cal{E}}^{\phi}(b(x),x_{2})c(x)\nonumber\\
&=&\Res_{x_{0}}x_{0}^{n} f(\phi(x,x_{0}+x_{2}),x)f(\phi(x,x_{0}),x)
Y_{\cal{E}}^{\phi}(Y_{\cal{E}}^{\phi}(a(x),x_{0})b(x), x_{2})c(x).\
\ \ \
\end{eqnarray*}
Since {\em $f(x,y)$ does not depend on $n$ and since $n$ is
arbitrary}, we have
\begin{eqnarray*}
& &f(\phi(x,x_{0}+x_{2}),x)f(\phi(x,x_{0}),x)
Y_{\cal{E}}^{\phi}(a(x),x_{0}+x_{2})Y_{\cal{E}}^{\phi}(b(x),x_{2})c(x)\nonumber\\
&=&f(\phi(x,x_{0}+x_{2}),x)f(\phi(x,x_{0}),x)
Y_{\cal{E}}^{\phi}(Y_{\cal{E}}^{e}(a(x),x_{0})b(x), x_{2})c(x).
\end{eqnarray*}
In view of Remark \ref{rcancelation}, we can multiply both sides by
$f(\phi(x,x_{0}),x)^{-1}$ to get
\begin{eqnarray}\label{enearfinal-new}
& &f(\phi(x,x_{0}+x_{2}),x)
Y_{\cal{E}}^{\phi}(a(x),x_{0}+x_{2})Y_{\cal{E}}^{\phi}(b(x),x_{2})c(x)\nonumber\\
&=&f(\phi(x,x_{0}+x_{2}),x)
Y_{\cal{E}}^{\phi}(Y_{\cal{E}}^{\phi}(a(x),x_{0})b(x), x_{2})c(x).
\end{eqnarray}
Write $f(\phi(x,z),x)=z^{k}g(x,z)$ for some $k\in \N,\; g(x,z)\in
\C((x))[[z]]$ with $g(x,0)\ne 0$. Then
\begin{eqnarray*}
f(\phi(x,x_{0}+x_{2}),x)=(x_{0}+x_{2})^{k}g(x,x_{0}+x_{2})
\end{eqnarray*}
and $g(x,x_{0}+x_{2})$ is a unit in $\C((x))[[x_{0},x_{2}]]$. By
cancellation we obtain
\begin{eqnarray*}
& &(x_{0}+x_{2})^{k}(Y_{\cal{E}}^{\phi}(a(x),x_{0}+x_{2})
Y_{\cal{E}}^{\phi}(b(x),x_{2})c(x)\\
&=&(x_{0}+x_{2})^{k}
Y_{\cal{E}}^{\phi}(Y_{\cal{E}}^{\phi}(a(x),x_{0})b(x),x_{2})c(x),
\end{eqnarray*}
as desired.

With $Y_{W}(a(x),z)=a(z)$ for $a(x)\in V$, we have
$Y_{W}(1_{W},z)=1_{W}$. Furthermore, for $a(x),b(x)\in V$, there
exists $h(x,y)\in \C[[x,y]]$ with  $h(\phi(x,z),x)\ne 0$ such that
$$h(x_{1},x_{2})a(x_{1})b(x_{2})\in \Hom (W,W((x_{1},x_{2}))).$$
Then
$$h(x_{1},x_{2})Y_{W}(a(x),x_{1})Y_{W}(b(x),x_{2})
=h(x_{1},x_{2})a(x_{1})b(x_{2})\in \Hom (W,W((x_{1},x_{2})))$$ and
$$h(\phi(x_{2},x_{0}),x_{2})(Y_{\cal{E}}^{\phi}(a(x),x_{0})b(x))|_{x=x_{2}}=
\left(h(x_{1},x_{2})a(x_{1})b(x_{2})\right)|_{x_{1}=\phi(x_{2},x_{0})}.$$
That is,
\begin{eqnarray*}
&&h(\phi(x_{2},x_{0}),x_{2})Y_{W}(Y_{\cal{E}}^{\phi}(a(x),x_{0})b(x),x_{2})
\nonumber\\
&=&\left(h(x_{1},x_{2})Y_{W}(a(x),x_{1})Y_{W}(b(x),x_{2})\right)|_{x_{1}=\phi(x_{2},x_{0})}.
\end{eqnarray*}
Therefore, $W$ is a $\phi$-coordinated quasi $V$-module. The
furthermore assertion is clear from the proof.
\end{proof}

Next, we are going to prove that every $\phi$-quasi compatible
subset of $\E(W)$ generates a nonlocal vertex algebra. To achieve
this goal, we first establish the following key result:

\bp{pgeneratingcomplicatedone} Let $\psi_{1}(x),\dots,\psi_{r}(x),
a(x),b(x),\phi_{1}(x), \dots,\phi_{s}(x) \in \E(W)$. Assume that the
ordered sequences $(a(x), b(x))$ and
$$(\psi_{1}(x),\dots,\psi_{r}(x), a(x),b(x),\phi_{1}(x),
\dots,\phi_{s}(x))$$
 are $\phi$-quasi compatible.
Then for any $n\in \Z$, the ordered sequence
$$(\psi_{1}(x),\dots,\psi_{r}(x),a(x)_{n}^{\phi}b(x),\phi_{1}(x),\dots,\phi_{s}(x))$$
is $\phi$-quasi compatible. The same assertion holds without the
prefix ``quasi.'' \ep

\begin{proof}
Let $f(x,y)\in \C[[x,y]]$ be such that $f(\phi(x,z),x)\ne 0$,
$$f(x_{1},x_{2})a(x_{1})b(x_{2})\in \Hom (W,W((x_{1},x_{2})))$$
and
\begin{eqnarray}\label{elong-exp}
& &\left(\prod_{1\le i<j\le r}f(y_{i},y_{j})\right)
\left(\prod_{1\le i\le r, 1\le j\le s}f(y_{i},z_{j})\right)
\left(\prod_{1\le i<j\le s}f(z_{i},z_{j})\right)\nonumber\\
& &\;\;\cdot f(x_{1},x_{2})
\left(\prod_{i=1}^{r}f(y_{i},x_{1})f(x_{2},y_{i})\right)
\left(\prod_{i=1}^{s}f(x_{1},z_{i})f(x_{2},z_{i})\right)\nonumber\\
& &\;\;\cdot \psi_{1}(y_{1})\cdots \psi_{r}(y_{r})
a(x_{1})b(x_{2})\phi_{1}(z_{1})\cdots \phi_{s}(z_{s})\nonumber\\
& &\in \Hom (W,W((y_{1},\dots,
y_{r},x_{1},x_{2},z_{1},\dots,z_{s}))).
\end{eqnarray}
Set
$$P=\prod_{1\le i<j\le r}f(y_{i},y_{j}),\;\;\;\;
Q=\prod_{1\le i<j\le s}f(z_{i},z_{j}),\;\;\;\; R=\prod_{1\le i\le
r,\; 1\le j\le s}f(y_{i},z_{j}).$$ {}From Proposition
\ref{pall-assoc} we have
$$\phi(x,z)=e^{zp(x)(d/dx)}x$$
for some $p(x)\in \C((x))$.
 Let $n\in
\Z$ be {\em arbitrarily fixed}. There exists a nonnegative integer
$k$ such that
\begin{eqnarray}\label{etruncationpsiphi}
x_{0}^{k+n}f(\phi(x,x_{0}),x)^{-1} \in \C((x))[[x_{0}]].
\end{eqnarray}
Using this and the fact
\begin{eqnarray}
\phi(\phi(x,z),-z)=\phi(x,0)=x=\phi(\phi(x,-z),z),
\end{eqnarray} we obtain
\begin{eqnarray*}\label{ecompatibilitythreeproof}
& &\prod_{i=1}^{r}f(x_{2},y_{i})^{k}
\prod_{j=1}^{s}f(x_{2},z_{j})^{k}
\psi_{1}(y_{1})\cdots\psi_{r}(y_{r})
(a(x)_{n}^{\phi}b(x))(x_{2})\phi_{1}(z_{1})\cdots\phi_{s}(z_{s})\nonumber\\
&=&\Res_{x_{0}}x_{0}^{n}\prod_{i=1}^{r}f(x_{2},y_{i})^{k}
\prod_{j=1}^{s}f(x_{2},z_{j})^{k}\nonumber\\
& &\ \ \ \ \cdot \psi_{1}(y_{1})\cdots \psi_{r}(y_{r})
(Y_{\cal{E}}^{\phi}(a(x_{2}),x_{0})b(x_{2}))\phi_{1}(z_{1})\cdots
\phi_{s}(z_{s})
\nonumber\\
&=&\Res_{x_{1}}\Res_{x_{0}}x_{0}^{n}\prod_{i=1}^{r}f(x_{2},y_{i})^{k}
\prod_{j=1}^{s}f(x_{2},z_{j})^{k}f(\phi(x_{2},x_{0}),x_{2})^{-1}
\nonumber\\
& &\ \ \ \ \cdot
x_{1}^{-1}\delta\left(\frac{\phi(x_{2},x_{0})}{x_{1}}\right)
\left(f(x_{1},x_{2})\psi_{1}(y_{1})\cdots \psi_{r}(y_{r})
a(x_{1})b(x_{2})\phi_{1}(z_{1})\cdots
\phi_{s}(z_{s})\right)\nonumber\\
&=&\Res_{x_{1}}\Res_{x_{0}}x_{0}^{n}\prod_{i=1}^{r}f(\phi(x_{1},-x_{0}),y_{i})^{k}
\prod_{j=1}^{s}f(\phi(x_{1},-x_{0}),z_{j})^{k}
f(\phi(x_{2},x_{0}),x_{2})^{-1}\nonumber\\
& &\ \ \ \ \cdot
x_{1}^{-1}\delta\left(\frac{\phi(x_{2},x_{0})}{x_{1}}\right)
\left(f(x_{1},x_{2})\psi_{1}(y_{1})\cdots \psi_{r}(y_{r})
a(x_{1})b(x_{2})\phi_{1}(z_{1})\cdots
\phi_{s}(z_{s})\right)\nonumber\\
&=&\Res_{x_{1}}\Res_{x_{0}}x_{0}^{n}
e^{-x_{0}p(x_{1})\frac{\partial}{\partial x_{1}}}
\left(\prod_{i=1}^{r}f(x_{1},y_{i})
\prod_{j=1}^{s}f(x_{1},z_{j})\right)^{k}f(\phi(x_{2},x_{0}),x_{2})^{-1}\nonumber\\
& &\ \ \ \ \cdot
x_{1}^{-1}\delta\left(\frac{\phi(x_{2},x_{0})}{x_{1}}\right)
\left(f(x_{1},x_{2})\psi_{1}(y_{1})\cdots \psi_{r}(y_{r})
a(x_{1})b(x_{2})\phi_{1}(z_{1})\cdots
\phi_{s}(z_{s})\right)\nonumber\\
&=&\Res_{x_{1}}\Res_{x_{0}}\sum_{t=0}^{k-1}\frac{(-1)^{t}}{t!}x_{0}^{n+t}
\left(p(x_{1})\frac{\partial}{\partial x_{1}}\right)^{t}
\left(\prod_{i=1}^{r}f(x_{1},y_{i})
\prod_{j=1}^{s}f(x_{1},z_{j})\right)^{k}\nonumber\\
& &\hspace{1cm}\cdot f(\phi(x_{2},x_{0}),x_{2})^{-1}
x_{1}^{-1}\delta\left(\frac{\phi(x_{2},x_{0})}{x_{1}}\right)\nonumber\\
& &\hspace{1cm}\cdot \left(f(x_{1},x_{2})\psi_{1}(y_{1})\cdots
\psi_{r}(y_{r}) a(x_{1})b(x_{2})\phi_{1}(z_{1})\cdots
\phi_{s}(z_{s})\right).
\end{eqnarray*}
Notice that for any power series $B$ and for $0\le t\le k-1$,
$\left(p(x_{1})\frac{\partial}{\partial x_{1}}\right)^{t}B^{k}$ is a
multiple of $B$. Using (\ref{elong-exp}) we have
\begin{eqnarray*}
& &P Q R \prod_{i=1}^{r}f(x_{2},y_{i})
\prod_{j=1}^{s}f(x_{2},z_{j})\\
&& \cdot \sum_{t=0}^{k-1}\frac{(-1)^{t}}{t!}x_{0}^{n+t}
\left(p(x_{1})\frac{\partial}{\partial x_{1}}\right)^{t}
\left(\prod_{i=1}^{r}f(y_{i},x_{1})
\prod_{j=1}^{s}f(x_{1},z_{j})\right)^{k}
f(\phi(x_{2},x_{0}),x_{2})^{-1}\nonumber\\
& &\cdot
x_{1}^{-1}\delta\left(\frac{\phi(x_{2},x_{0})}{x_{1}}\right)
\left(f(x_{1},x_{2})\psi_{1}(y_{1})\cdots \psi_{r}(y_{r})
a(x_{1})b(x_{2})\phi_{1}(z_{1})\cdots
\phi_{s}(z_{s})\right)\nonumber\\
&\in& \left(\Hom (W,W((y_{1},\dots,
y_{r},x_{2},z_{1},\dots,z_{s}))\right)
((x_{0}))[[x_{1},x_{1}^{-1}]].
\end{eqnarray*}
Then
\begin{eqnarray}
& &P Q R\prod_{i=1}^{r}f(x_{2},y_{i})^{k+1}
\prod_{j=1}^{s}f(x_{2},z_{j})^{k+1}\nonumber\\
&&\ \ \ \ \cdot \psi_{1}(y_{1})\cdots\psi_{r}(y_{r})
(a(x)_{n}^{\phi}b(x))(x_{2})\phi_{1}(z_{1})\cdots\phi_{s}(z_{s})\nonumber\\
&\in& \Hom (W,W((y_{1},\dots, y_{r},x_{2},z_{1},\dots,z_{s}))).
\end{eqnarray}
This proves that the sequence
$(\psi_{1}(x),\dots,\psi_{r}(x),a(x)_{n}^{\phi}b(x), \phi_{1}(x),
\dots,\phi_{s}(x))$ is $\phi$-quasi compatible. The last assertion
also follows from the proof.
\end{proof}

The following is the main result of this section:

\bt{tfirst} Let $W$ be a vector space, $\phi(x,z)$ an associate of
the additive formal group $F_{\rm a}(x,y)$, and $U$ a (resp.
$\phi$-quasi) compatible subset of $\E(W)$. There exists a
$Y_{\E}^{\phi}$-closed (resp. $\phi$-quasi) compatible subspace of
$\E(W)$, that contains $U$ and $1_{W}$. Denote by $\<U\>_{\phi}$ the
smallest such subspace. Then $(\<U\>_{\phi},Y_{\E}^{\phi},1_{W})$
carries the structure of a nonlocal vertex algebra and $W$ is a
$\phi$-coordinated (resp. quasi) $\<U\>_{\phi}$-module with
$Y_{W}(\alpha(x),z)=\alpha(z)$ for $\alpha(x)\in \<U\>_{\phi}$.
 \et

\begin{proof} By Zorn's lemma, there exists a maximal
quasi compatible subspace $V$ of $\E(W)$, containing both $U$ and
$1_{W}$. It follows from Proposition \ref{pgeneratingcomplicatedone}
that $V$ is $Y_{\E}^{\phi}$-closed. This proves the first assertion.
Furthermore, by Theorem \ref{tclosed}, $(V,Y_{\E}^{\phi},1_{W})$
carries the structure of a nonlocal vertex algebra with $W$ as a
$\phi$-coordinated quasi module. By definition, $\<U\>_{\phi}$ is
the intersection of all $Y_{\E}^{\phi}$-closed  (resp. $\phi$-quasi)
compatible subspaces of $\E(W)$, containing  both $U$ and $1_{W}$.
The rest follows from Theorem \ref{tclosed}.
\end{proof}

Just as with usual quasi modules for a nonlocal vertex algebra, the
state-field correspondence for $\phi$-coordinated quasi modules is
also a homomorphism.

\bp{pstate-field} Let $V$ be a nonlocal vertex algebra and let
$(W,Y_{W})$ be a $\phi$-coordinated quasi $V$-module. Then
\begin{eqnarray}
Y_{W}(Y(u,x_{0})v,x)=Y_{\E}^{\phi}(u(x),x_{0})v(x)
\end{eqnarray}
for $u,v\in V$, where $u(x)=Y_{W}(u,x),\; v(x)=Y_{W}(v,x)\in \E(W)$.
\ep

\begin{proof} For $u,v\in V$, there exists
$p(x,y)\in \C[[x,y]]$ such that $p(\phi(x,z),x)\ne 0$,
\begin{eqnarray}\label{epuvcondition}
p(x_{1},x_{2})u(x_{1})v(x_{2})\in \Hom(W,W((x_{1},x_{2})))
\end{eqnarray}
 and
$$p(\phi(x,x_{0}),x)Y_{W}(Y(u,x_{0})v,x)
=\left(p(x_{1},x)u(x_{1})v(x)\right)|_{x_{1}=\phi(x,x_{0})}.$$ With
(\ref{epuvcondition}), we also have
$$p(\phi(x,x_{0}),x)Y_{\E}^{\phi}(u(x),x_{0})v(x)
=\left(p(x_{1},x)u(x_{1})v(x)\right)|_{x_{1}=\phi(x,x_{0})}.$$ Thus
$$p(\phi(x,x_{0}),x)Y_{W}(Y(u,x_{0})v,x)
=p(\phi(x,x_{0}),x)Y_{\E}^{\phi}(u(x),x_{0})v(x).$$ As the powers of
$x_{0}$ in both $Y_{W}(Y(u,x_{0})v,x)$ and
$Y_{\E}^{\phi}(u(x),x_{0})v(x)$ are lower truncated, with
$p(\phi(x,x_{0}),x)\in \C((x))[[x_{0}]]$ nonzero we obtain the
desired relation by cancellation.
\end{proof}

\br{roldcase} {\em Consider the case with $\phi(x,z)=x+z$. {}From
Lemma \ref{lphi-quasicompatible}, $\phi$-quasi compatibility is the
same as quasi compatibility. Furthermore, for a quasi compatible
pair $(a(x),b(x))$ in $\E(W)$, we have
$Y_{\E}^{\phi}(a(x),z)b(x)=Y_{\E}(a(x),z)b(x)$, which was defined in
\cite{li-qva1} by
$$Y_{\E}(a(x),z)b(x)
=\iota_{x,z}(1/p(x+z,x))\left(p(x_{1},x)a(x_{1})b(x)\right)|_{x_{1}=x+z},$$
where $p(x,y)$ is any nonzero element of $\C[[x,y]]$ such that
$$p(x_{1},x_{2})a(x_{1})b(x_{2})\in \Hom (W,W((x_{1},x_{2}))).$$
On the other hand, with $\phi(x,z)=x+z$, a $\phi$-coordinated quasi
module a nonlocal vertex algebra is simply a quasi module. In view
of these, Theorem \ref{tfirst} generalizes the corresponding results
of \cite{li-qva1}. } \er

\br{rtrivialcase} {\em Consider the extreme case with $\phi(x,z)=x$.
A pair $(a(x),b(x))$ in $\E(W)$ is $\phi$-quasi compatible if and
only if there exists $p(x_{1},x_{2})\in \C[[x_{1},x_{2}]]$ with
$p(x_{2},x_{2})\ne 0$ such that
$$p(x_{1},x_{2})a(x_{1})b(x_{2})\in \Hom (W,W((x_{1},x_{2}))).$$
Assuming that $(a(x),b(x))$ is $\phi$-quasi compatible with
$p(x_{1},x_{2})\in \C[[x_{1},x_{2}]]$ satisfying the above
condition, we have
$$Y_{\E}^{\phi}(a(x),z)b(x)=p(x,x)^{-1}\left(p(x_{1},x)a(x_{1})b(x)\right)|_{x_{1}=x},$$
which is independent of $z$, where $p(x,x)^{-1}$ stands for the
inverse of $p(x,x)$ in $\C((x))$. Then  the nonlocal vertex algebra
$\<U\>_{\phi}$, associated to a $\phi$-quasi compatible subset $U$
of $\E(W)$ by Theorem \ref{tfirst}, is in fact an ordinary
associative algebra. {}From Theorem \ref{tfirst}, the vector space
$W$ is a $\phi$-coordinated quasi module for $\<U\>_{\phi}$ viewed
as a nonlocal vertex algebra, but $W$ is not a module in the usual
sense for $\<U\>_{\phi}$ viewed as either an associative algebra or
a nonlocal vertex algebra in general.} \er

\section{$\phi$-coordinated modules for weak quantum vertex
algebras with $\phi(x,z)=xe^{z}$}

In this section, we study $\phi$-coordinated quasi modules for weak
quantum vertex algebras with $\phi$ specialized to
$\phi(x,z)=xe^{z}$. We first continue with Section 4 to formulate
notions of (quasi) $\S_{trig}$-local subset and prove that the
nonlocal vertex algebra generated by any quasi $\S_{trig}$-local
subset is a weak quantum vertex algebra. We then present certain
axiomatic results on $\phi$-coordinated modules for weak quantum
vertex algebras. In particular, we establish a Jacobi-type identity.

Let $W$ be a vector space as in Section 4. Throughout this section,
we assume $\phi(x,z)=xe^{z}$ and we denote $Y_{\E}^{\phi}$ by
$Y_{\E}^{e}$. That is,
\begin{eqnarray}
Y_{\E}^{e}(a(x),z)b(x)
=p(xe^{z},x)^{-1}\left(p(x_{1},x)a(x_{1})b(x)\right)|_{x_{1}=xe^{z}}
\end{eqnarray}
for any quasi compatible pair $(a(x),b(x))$ in $\E(W)$ with nonzero
$p(x_{1},x_{2})\in \C[[x_{1},x_{2}]]$ such that
(\ref{epab-condition-phi}) holds.

First we formulate the following notions:

\bd{dsemi-local} {\em  A subset $U$ of $\E(W)$ is said to be {\em
$\S_{trig}$-local} if for any $a(x),b(x)\in U$, there exist
$$ u_{i}(x),v_{i}(x)\in U,\ q_{i}(x)\in \C(x)\ \ (i=1,\dots,r),$$
where $\C(x)$ denotes the field of rational functions, such that
\begin{eqnarray}\label{etrig-slocal}
&&(x_{1}-x_{2})^{k}a(x_{1})b(x_{2})
=(x_{1}-x_{2})^{k}\sum_{i=1}^{r}\iota_{x_{2},x_{1}}(q_{i}(x_{1}/x_{2}))
u_{i}(x_{2})v_{i}(x_{1})
\end{eqnarray}
for some $k\in \N$. The notion of {\em quasi $\S_{trig}$-local
subset} is defined by weakening the above condition as
\begin{eqnarray}\label{esemi-pab}
&&p(x_{1}/x_{2})a(x_{1})b(x_{2})
=p(x_{1}/x_{2})\sum_{i=1}^{r}\iota_{x_{2},x_{1}}(q_{i}(x_{1}/x_{2}))
u_{i}(x_{2})v_{i}(x_{1})
\end{eqnarray}
for some nonzero polynomial $p(x)\in \C[x_{1},x_{2}]$. } \ed

These notions single out a family of compatible subsets and a family
of quasi compatible subsets as we show next.

\bl{lpracticalcase} Every (resp. quasi) $\S_{trig}$-local subset of
$\E(W)$ is (resp. quasi) compatible. \el

\begin{proof} Let us first consider the quasi case. Let $U$ be a quasi $\S_{trig}$-local subset.
We must prove that every finite sequence in $U$ is quasi compatible.
To prove this we use induction on the length $n$ of sequences. Let
$(a(x),b(x))$ be an ordered pair in $U$. By assumption, there exist
$0\ne p(x)\in \C[x]$, $a^{(i)}(x), b^{(i)}(x)\in U$ and $q_{i}(x)\in
\C(x)$ for $i=1,\dots,r$ such that
\begin{eqnarray}\label{e4.3}
p(x_{1}/x_{2})a(x_{1})b(x_{2})=p(x_{1}/x_{2})
\sum_{i=1}^{r}\iota_{x_{2},x_{1}}(q_{i}(x_{1}/x_{2}))b^{(i)}(x_{2})a^{(i)}(x_{1}).
\end{eqnarray}
The expression on the left-hand side lies in $\Hom
(W,W((x_{1}))((x_{2})))$ while the expression on the right-hand side
lies in $\Hom (W,W((x_{2}))((x_{1})))$. This forces the expressions
on both sides to lie in $\Hom (W,W((x_{1},x_{2})))$. Thus
$(a(x),b(x))$ is quasi compatible, proving the case for $n=2$.

Now assume that $n\ge 2$ and that any sequence  in $U$ of length $n$
is quasi compatible. Let $\psi^{(1)}(x),\dots,\psi^{(n+1)}(x)\in U$.
From the inductive hypothesis, there exists  $0\ne f(x)\in \C[x]$
such that
\begin{eqnarray}\label{eproduct1}
\left(\prod_{2\le i<j\le n+1}f(x_{i}/x_{j})\right)
\psi^{(2)}(x_{2})\cdots \psi^{(n+1)}(x_{n+1}) \in \Hom
(W,W((x_{2},\dots,x_{n+2}))).
\end{eqnarray}
By assumption there exist $0\ne p(x)\in \C[x]$, $a^{(i)}(x),
b^{(i)}(x)\in U$ and $q_{i}(x)\in \C(x)$ for $i=1,\dots,r$ such that
\begin{eqnarray}\label{eproduct2}
p(x_{1}/x_{2})\psi^{(1)}(x_{1})\psi^{(2)}(x_{2})
=p(x_{1}/x_{2})\sum_{i=1}^{r}\iota_{x_{2},x_{1}}(q_{i}(x_{1}/x_{2}))b^{(i)}(x_{2})a^{(i)}(x_{1}).
\end{eqnarray}
{}From the inductive hypothesis again, there exists $0\ne g(x)\in
\C[x]$ such that
\begin{eqnarray}\label{eproduct3}
& &\left( \prod_{1\le i<j\le n+1,\; i,j\ne 2}g(x_{i}/x_{j})\right)
a^{(s)}(x_{1})\psi^{(3)}(x_{3})\cdots \psi^{(n+1)}(x_{n+1})
\nonumber\\
&\in& \Hom (W, W((x_{1},x_{3},x_{4},\dots,x_{n+1})))
\end{eqnarray}
for $s=1,\dots,r$. Using (\ref{eproduct2}) we have
\begin{eqnarray}\label{eproduct4}
& &\left(\prod_{2\le i<j\le n+1}f(x_{i}/x_{j}) \prod_{1\le i<j\le
n+1,\; i,j\ne 2}g(x_{i}/x_{j})\right)
p(x_{1}/x_{2})\psi^{(1)}(x_{1})\cdots \psi^{(n+1)}(x_{n+1})\nonumber\\
&=&\left(\prod_{2\le i<j\le n+1}f(x_{i}/x_{j})
\prod_{1\le i<j\le n+1,\; i,j\ne 2}g(x_{i}/x_{j})\right)
p(x_{1}/x_{2})\nonumber\\
& &\ \ \ \ \ \cdot
\sum_{s=1}^{r}\iota_{x_{2},x_{1}}(q_{i}(x_{1}/x_{2}))
b^{(s)}(x_{2})a^{(s)}(x_{1})\psi^{(3)}(x_{3})\cdots
\psi^{(n+1)}(x_{n+1}).
\end{eqnarray}
{}From (\ref{eproduct1}), the expression on the left-hand side of
(\ref{eproduct4}) lies in
$$\Hom (W,W((x_{1}))((x_{2},x_{3},x_{4},\dots,x_{n+1}))),$$
and by (\ref{eproduct3}), the expression on the right-hand side of
(\ref{eproduct4}) lies in
$$\Hom (W,W((x_{2}))((x_{1},x_{3},x_{4},\dots,x_{n+1}))).$$
This forces the expressions on both sides to lie in the space $$\Hom
(W,W((x_{1},x_{2},x_{3},x_{4},\dots,x_{n+1}))).$$ In particular,
the expression on the left-hand side (\ref{eproduct4}) lies in
$$\Hom (W,W((x_{1},x_{2},x_{3},x_{4},\dots,x_{n+1}))).$$ This proves
that the sequence $(\psi^{(1)}(x),\dots,\psi^{(n+1)}(x))$ is quasi
compatible, completing the induction. {}From the proof, it is clear
 that if $U$ is $\S_{trig}$-local, then $U$ is compatible.
\end{proof}

In view of Lemma \ref{lpracticalcase} and Theorem \ref{tfirst}, for
any quasi $\S_{trig}$-local subset $U$ of $\E(W)$, we have a
nonlocal vertex algebra $\<U\>_{\phi}$ generated by $U$ with
$\phi(x,z)=xe^{z}$. In the following we are going to prove that
$\<U\>_{\phi}$ is a weak quantum vertex algebra. To this end we
first prove:

\bp{pconvert} Let $V$ be a $Y_{\E}^{e}$-closed quasi compatible
subspace of $\E(W)$. Suppose
$$a(x),b(x),u_{i}(x),v_{i}(x)\in V,\; 0\ne p(x)\in \C[x],\; q_{i}(x)\in
\C(x)\ \ (i=1,\dots,r)$$ satisfy
\begin{eqnarray}\label{econdition}
p(x_{1}/x_{2})a(x_{1})b(x_{2})
=\sum_{i=1}^{r}p(x_{1}/x_{2})\iota_{x_{2},x_{1}}(q_{i}(x_{1}/x_{2}))
u_{i}(x_{2})v_{i}(x_{1}).
\end{eqnarray}
 Then
\begin{eqnarray}\label{epfsum01}
& &p(e^{x_{1}-x_{2}})Y_{\E}^{e}(a(x),x_{1})Y_{\E}^{e}(b(x),x_{2})
\nonumber\\
&=&p(e^{x_{1}-x_{2}})\sum_{i=1}^{r}\iota_{x_{2},x_{1}}(q_{i}(e^{x_{1}-x_{2}}))
Y_{\E}^{e}(u_{i}(x),x_{2})Y_{\E}^{e}(v_{i}(x),x_{1}).
\end{eqnarray}
Furthermore, we have
\begin{eqnarray}\label{epfsum}
& &(x_{1}-x_{2})^{k}Y_{\E}^{e}(a(x),x_{1})Y_{\E}^{e}(b(x),x_{2})
\nonumber\\
&=&(x_{1}-x_{2})^{k}\sum_{i=1}^{r}\iota_{x_{2},x_{1}}(q_{i}(e^{x_{1}-x_{2}}))
Y_{\E}^{e}(u_{i}(x),x_{2})Y_{\E}^{e}(v_{i}(x),x_{1}),
\end{eqnarray}
where $k$ is the multiplicity of the zero of $p(x)$ at $x=1$. \ep

\begin{proof} Let $\theta(x)\in V$ be arbitrarily fixed.
There exists $0\ne f(x_{1},x_{2})\in \C[[x_{1},x_{2}]]$ such that
\begin{eqnarray*}
& &f(z,x)b(z)\theta(x)\in \Hom (W,W((x,z))),\\
& &f(y,z)f(y,x)f(z,x)a(y)b(z)\theta(x)\in \Hom (W,W((x,y,z))).
\end{eqnarray*}
By Lemma \ref{lclosed}, we have
\begin{eqnarray*}
& &f(xe^{x_{1}},xe^{x_{2}})f(xe^{x_{1}},x)f(xe^{x_{2}},x)
Y_{\E}^{e}(a(x),x_{1})Y_{\E}^{e}(b(x),x_{2})\theta (x)\\
&=&\left(f(y,z)f(y,x)f(z,x) a(y)b(z)\theta
(x)\right)|_{y=xe^{x_{1}},z=xe^{x_{2}}}.
\end{eqnarray*}
Replacing $f(x_{1},x_{2})$ with a multiple of $f(x_{1},x_{2})$ if
necessary, we also have
\begin{eqnarray*}
& &f(xe^{x_{1}},xe^{x_{2}})f(xe^{x_{1}},x)f(xe^{x_{2}},x)
Y_{\E}^{e}(u_{i}(x),x_{2})Y_{\E}^{e}(v_{i}(x),x_{1})\theta (x)
\nonumber\\
&=&\left(f(y,z)f(y,x)f(z,x) u_{i}(z)v_{i}(y)\theta
(x)\right)|_{z=xe^{x_{2}},y=xe^{x_{1}}}
\end{eqnarray*}
for $i=1,\dots,r$. Therefore,
\begin{eqnarray*}
&&f(xe^{x_{1}},xe^{x_{2}})f(xe^{x_{1}},x)f(xe^{x_{2}},x)
p(e^{x_{1}-x_{2}})
Y_{\E}^{e}(a(x),x_{1})Y_{\E}^{e}(b(x),x_{2})\theta (x)\nonumber\\
&=&\left(f(y,z)f(y,x)f(z,x)p(y/z) a(y)b(z)\theta
(x)\right)|_{y=xe^{x_{1}},z=xe^{x_{2}}}\\
&=&\sum_{i=1}^{r}\left(f(y,z)f(y,x)f(z,x)p(y/z)q_{i}(y/z)
u_{i}(z)v_{i}(y)\theta
(x)\right)|_{z=xe^{x_{2}},y=xe^{x_{1}}}\\
&=&f(xe^{x_{1}},xe^{x_{2}})f(xe^{x_{1}},x)f(xe^{x_{2}},x)p(e^{x_{1}-x_{2}})\nonumber\\
& &\ \ \ \ \cdot
\sum_{i=1}^{r}\iota_{x_{2},x_{1}}(q_{i}(e^{x_{1}-x_{2}}))Y_{\E}^{e}(u_{i}(x),x_{2})
Y_{\E}^{e}(v_{i}(x),x_{1})\theta (x).
\end{eqnarray*}
In view of Remark \ref{rcancelation}, (by cancellation) we have
\begin{eqnarray*}
&&f(xe^{x_{1}},xe^{x_{2}}) p(e^{x_{1}-x_{2}})
Y_{\E}^{e}(a(x),x_{1})Y_{\E}^{e}(b(x),x_{2})\theta (x)\nonumber\\
&=&f(xe^{x_{1}},xe^{x_{2}})p(e^{x_{1}-x_{2}})
\sum_{i=1}^{r}\iota_{x_{2},x_{1}}(q_{i}(e^{x_{1}-x_{2}}))Y_{\E}^{e}(u_{i}(x),x_{2})
Y_{\E}^{e}(v_{i}(x),x_{1})\theta (x).\ \ \ \
\end{eqnarray*}
Write
$f(z_{1},z_{2})p(z_{1}/z_{2})=(z_{1}-z_{2})^{s}z_{2}^{-s'}g(z_{1},z_{2})$
with $s,s'\in \N,\; g(z_{1},z_{2})\in \C[[z_{1},z_{2}]]$ such that
$g(z,z)\ne 0$. Then
$$f(xe^{x_{1}},xe^{x_{2}})
p(e^{x_{1}-x_{2}})=(x_{1}-x_{2})^{s}x^{s}E(x_{1},x_{2})^{s}
(xe^{x_{2}})^{-s'}g(xe^{x_{1}},xe^{x_{2}}),$$ where
$E(x_{1},x_{2})=\sum_{n\ge 1}\frac{1}{n!}(x_{1}^{n-1}-x_{2}^{n-1})$
is a unit in $\C[[x_{1},x_{2}]]$.
 Noticing that $g(xe^{x_{1}},xe^{x_{2}})$ is a unit in
$\C((x))[[x_{1},x_{2}]]$, by cancellation we get
\begin{eqnarray*}
& &(x_{1}-x_{2})^{s}Y_{\E}^{e}(a(x),x_{1})Y_{\E}^{e}(b(x),x_{2})\nonumber\\
&=&\sum_{i=1}^{r}(x_{1}-x_{2})^{s}\iota_{x_{2},x_{1}}(q_{i}(e^{x_{1}-x_{2}}))
Y_{\E}^{e}(u_{i}(x),x_{2})Y_{\E}^{e}(v_{i}(x),x_{1}).
\end{eqnarray*}
Combining this with weak associativity (Theorem \ref{tclosed}) we
obtain
\begin{eqnarray}\label{e-jacobi-classical}
& &x_{0}^{-1}\delta\left(\frac{x_{1}-x_{2}}{x_{0}}\right)
Y_{\E}^{e}(a(x),x_{1})Y_{\E}^{e}(b(x),x_{2})\nonumber\\
& &\ \ \ \ -x_{0}^{-1}\delta\left(\frac{x_{2}-x_{1}}{-x_{0}}\right)
\sum_{i=1}^{r}\iota_{x_{2},x_{1}}(q_{i}(e^{x_{1}-x_{2}}))
Y_{\E}^{e}(u_{i}(x),x_{2})Y_{\E}^{e}(v_{i}(x),x_{1})\nonumber\\
&=&x_{2}^{-1}\delta\left(\frac{x_{1}-x_{0}}{x_{2}}\right)
Y_{\E}^{e}(Y_{\E}^{e}(a(x),x_{0})b(x),x_{2}).
\end{eqnarray}
{}From (\ref{econdition}) we have
$$p(x_{1}/x_{2})a(x_{1})b(x_{2})\in \Hom (W,W((x_{1},x_{2}))),$$
so that
$$p(e^{x_{0}})Y_{\E}^{e}(a(x),x_{0})b(x)
=\left(p(x_{1}/x)a(x_{1})b(x)\right)|_{x_{1}=xe^{x_{0}}},$$ which
involves only nonnegative integer powers of $x_{0}$. Multiplying the
both sides of (\ref{e-jacobi-classical}) by $p(e^{x_{0}})$ and then
taking $\Res_{x_{0}}$ we get
\begin{eqnarray*}
&&p(e^{x_{1}-x_{2}})Y_{\E}^{e}(a(x),x_{1})Y_{\E}^{e}(b(x),x_{2})
\nonumber\\
&=&p(e^{x_{1}-x_{2}})\sum_{i=1}^{r}\iota_{x_{2},x_{1}}(q_{i}(e^{x_{1}-x_{2}}))
Y_{\E}^{e}(u_{i}(x),x_{2})Y_{\E}^{e}(v_{i}(x),x_{1}).
\end{eqnarray*}
Let $k$ be the multiplicity of the zero of $p(x)$ at $x=1$. Then
$p(e^{z})=z^{k}h(z)$ where $h(z)\in \C[[z]]$ with $h(0)\ne 0$. By
cancellation we obtain
\begin{eqnarray*}
&&(x_{1}-x_{2})^{k}Y_{\E}^{e}(a(x),x_{1})Y_{\E}^{e}(b(x),x_{2})\nonumber\\
&=&\sum_{i=1}^{r}(x_{1}-x_{2})^{k}\iota_{x_{2},x_{1}}(q_{i}(e^{x_{1}-x_{2}}))
Y_{\E}^{e}(u_{i}(x),x_{2})Y_{\E}^{e}(v_{i}(x),x_{1}),
\end{eqnarray*}
as desired.
\end{proof}

As the main result of this section we have:

\bt{tqva} Let $W$ be a vector space and let $U$ be any (resp. quasi)
$\S_{trig}$-local subset of $\E(W)$. Then $\<U\>_{\phi}$ is a weak
quantum vertex algebra and $W$ is a $\phi$-coordinated (resp. quasi)
module with $\phi(x,z)=xe^{z}$. \et

\begin{proof} We only need to prove that $\<U\>_{\phi}$ is a weak quantum vertex algebra.
As $\<U\>_{\phi}$ is the smallest $Y_{\E}^{e}$-closed quasi
compatible subspace containing $U$ and $1_{W}$, we see that
$\<U\>_{\phi}$ as a nonlocal vertex algebra is generated by $U$.
Given that $U$ is quasi $\S_{trig}$-local, from Proposition
\ref{pconvert}, we have that
$$\{Y_{\E}^{e}(a(x),z)\;|\; a(x)\in U\}$$
is an $\S$-local subset of $\E(\<U\>_{\phi})$ in the sense of
\cite{li-qva1}. Then by \cite{ltw} (Proposition 2.6), $\<U\>_{\phi}$
is a weak quantum vertex algebra.
\end{proof}

\br{rqaalgebra} {\em Let $W$ be a highest weight module for a
quantum affine algebra $U_{q}(\hat{\g})$ with $q$ a complex number
(see \cite{drinfeld}, \cite{fj}). It is straightforward to see that
the generating functions of the generators in the Drinfeld
realization form a quasi $\S_{trig}$-local subset $U_{W}$ of
$\E(W)$. By Theorem \ref{tqva}, $U_{W}$ generates a weak quantum
vertex algebra with $W$ as a $\phi$-coordinated quasi module where
$\phi(x,z)=xe^{z}$. In a sequel, we shall study the associated weak
quantum vertex algebras in detail. } \er

We next study $\phi$-coordinated quasi modules for a general weak
quantum vertex algebra.

\bp{palgebra-module} Let $V$ be a nonlocal vertex algebra and let
$(W,Y_{W})$ be a $\phi$-coordinated quasi $V$-module. Assume that
$$u,v, u^{(i)},v^{(i)}\in V,\; f_{i}(x)\in \C(x) \ \ (i=1,\dots,r)$$
satisfy the relation
\begin{eqnarray}\label{eyab=uv}
&&(x_{1}-x_{2})^{k}Y(u,x_{1})Y(v,x_{2})\nonumber\\
&=&(x_{1}-x_{2})^{k}\sum_{i=1}^{r}\iota_{x_{2},x_{1}}(f_{i}(e^{x_{1}-x_{2}}))
Y(v^{(i)},x_{2})Y(u^{(i)},x_{1})
\end{eqnarray}
for some $k\in \N$. Suppose that $p(x_{1},x_{2})\in
\C[[x_{1},x_{2}]]$ is nonzero such that
$$p(x_{1},x_{2})Y_{W}(u,x_{1})Y_{W}(v,x_{2})\in \Hom (W,W((x_{1},x_{2}))).$$
Then
\begin{eqnarray}
&&p(x_{1},x_{2})Y_{W}(u,x_{1})Y_{W}(v,x_{2})\nonumber\\
&=&p(x_{1},x_{2})
\sum_{i=1}^{r}\iota_{x_{2},x_{1}}(f_{i}(x_{1}/x_{2}))Y_{W}(v^{(i)},x_{2})Y_{W}(u^{(i)},x_{1}).
\end{eqnarray}
 \ep

\begin{proof} With (\ref{eyab=uv}), by Corollary 5.3 of \cite{li-qva1} we have
$$Y(u,x)v=\sum_{i=1}^{r}\iota_{x,0}(f_{i}(e^{x}))e^{x\D}Y(v^{(i)},-x)u^{(i)}.$$
{}From definition, there exists a nonzero polynomial
$q(x_{1},x_{2})$ such that
$$q(x_{1},x_{2})Y_{W}(u,x_{1})Y_{W}(v,x_{2})\in \Hom
(W,W((x_{1},x_{2})))$$ and such that
$q(x_{1},x_{2})f_{i}(x_{1}/x_{2})\in \C[x_{1},x_{2}]$,
\begin{eqnarray*}
q(x_{1},x_{2})f_{i}(x_{1}/x_{2})Y_{W}(v^{(i)},x_{2})Y_{W}(u^{(i)},x_{1})\in
\Hom (W,W((x_{1},x_{2})))
\end{eqnarray*}
for $i=1,\dots,r$. Then, using Lemma \ref{lDproperty} we get
\begin{eqnarray*}
&&\left(q(x_{1},x_{2})Y_{W}(u,x_{1})Y_{W}(v,x_{2})\right)|_{x_{1}=x_{2}e^{x_{0}}}\\
&=&q(x_{2}e^{x_{0}},x_{2})Y_{W}(Y(u,x_{0})v,x_{2})\\
&=&\sum_{i=1}^{r}(q(x_{2}e^{x_{0}},x_{2})f_{i}(e^{x_{0}}))
Y_{W}(e^{x_{0}\D}Y(v^{(i)},-x_{0})u^{(i)},x_{2})\\
&=&\sum_{i=1}^{r}(q(x_{2}e^{x_{0}},x_{2})f_{i}(e^{x_{0}}))Y_{W}(Y(v^{(i)},-x_{0})u^{(i)},x_{2}e^{x_{0}}).
\end{eqnarray*}
We also have
\begin{eqnarray*}
&&\left(q(x_{1},x_{2})\sum_{i=1}^{r}f_{i}(x_{1}/x_{2})
Y_{W}(v^{(i)},x_{2})Y_{W}(u^{(i)},x_{1})\right)|_{x_{2}=x_{1}e^{-x_{0}}}\\
&=&\sum_{i=1}^{r}(q(x_{1},x_{1}e^{-x_{0}})f_{i}(e^{x_{0}}))Y_{W}(Y(v^{(i)},-x_{0})u^{(i)},x_{1}).
\end{eqnarray*}
Then using Remark \ref{rsimplefact} we have
\begin{eqnarray*}
&&\left(q(x_{1},x_{2})Y_{W}(u,x_{1})Y_{W}(v,x_{2})\right)|_{x_{1}=x_{2}e^{x_{0}}}\\
&=&\left(\left(q(x_{1},x_{2})\sum_{i=1}^{r}f_{i}(x_{1}/x_{2})
Y_{W}(v^{(i)},x_{2})Y_{W}(u^{(i)},x_{1})\right)|_{x_{2}=x_{1}e^{-x_{0}}}\right)|_{x_{1}=x_{2}e^{x_{0}}}
\\
&=&\left(q(x_{1},x_{2})\sum_{i=1}^{r}f_{i}(x_{1}/x_{2})
Y_{W}(v^{(i)},x_{2})Y_{W}(u^{(i)},x_{1})\right)|_{x_{1}=x_{2}e^{x_{0}}}.
\end{eqnarray*}
Using Remark \ref{rsimplefact} again we get
\begin{eqnarray*}
q(x_{1},x_{2})Y_{W}(u,x_{1})Y_{W}(v,x_{2})
=\sum_{i=1}^{r}q(x_{1},x_{2})f_{i}(x_{1}/x_{2})
Y_{W}(v^{(i)},x_{2})Y_{W}(u^{(i)},x_{1}).
\end{eqnarray*}
Then
\begin{eqnarray*}
&&q(x_{1},x_{2})\left(p(x_{1},x_{2})Y_{W}(u,x_{1})Y_{W}(v,x_{2})\right)\\
&=&q(x_{1},x_{2})\left(p(x_{1},x_{2})\sum_{i=1}^{r}\iota_{x_{2},x_{1}}(f_{i}(x_{1}/x_{2}))
Y_{W}(v^{(i)},x_{2})Y_{W}(u^{(i)},x_{1})\right).
\end{eqnarray*}
Multiplying both sides by the inverse of $q(x_{1},x_{2})$ in
$\C((x_{2}))((x_{1}))$  we obtain the desired relation.
\end{proof}

\br{rlog-exp} {\em Note that for any $f(x)\in x\C[[x]],\; g(x)\in
\C((x))$, the composition $g(f(x))$ exists in $\C((x))$. Set
\begin{eqnarray}
\log(1+x)=\sum_{n\ge 1}(-1)^{n-1}\frac{1}{n}x^{n}\in x\C[[x]].
\end{eqnarray}
For any $f(x)\in x\C[[x]]$, we have
$$\log(1+f(x))=\sum_{n\ge 1}(-1)^{n-1}\frac{1}{n}f(x)^{n}\in
x\C[[x]],$$
$$e^{f(x)}=\sum_{n\ge 0}\frac{1}{n!}f(x)^{n}\in \C[[x]].$$
Using formal calculus one can show
\begin{eqnarray}
e^{\log(1+z)}=1+z\ \ \mbox{ and }\ \ \log(1+(e^{x}-1))=x.
\end{eqnarray}
Let $E(x_{1},x_{2})\in \C((x_{1},x_{2}))$. Set
$$F(x_{0},x_{2})=E(x_{2}e^{x_{0}},x_{2})\in \C((x_{2}))[[x_{0}]].$$
Then
$$F(\log(1+z),x_{2})=E(x_{2}(1+z),x_{2})\in \C((x_{2}))[[z]].$$}
\er

\bl{ldelta-identity} Let $W$ be any vector space and let
\begin{eqnarray*}
&&A(x_{1},x_{2})\in \Hom(W,W((x_{1}))((x_{2}))),\ \
B(x_{1},x_{2})\in \Hom(W,W((x_{2}))((x_{1}))),\\
&&\hspace{2.5cm}C(x_{0},x_{2})\in (\Hom (W,W((x_{2}))))((x_{0})).
\end{eqnarray*}
If there exists a nonnegative integer $k$ such that
\begin{eqnarray*}
&&(x_{1}-x_{2})^{k}A(x_{1},x_{2})=(x_{1}-x_{2})^{k}B(x_{1},x_{2}),\\
&&\left((x_{1}-x_{2})^{k}A(x_{1},x_{2})\right)|_{x_{1}=x_{2}e^{x_{0}}}
=x_{2}^{k}(e^{x_{0}}-1)^{k}C(x_{0},x_{2}),
\end{eqnarray*}
then
\begin{eqnarray}\label{edelta-identity}
&&(x_{2}z)^{-1}\delta\left(\frac{x_{1}-x_{2}}{x_{2}z}\right)
A(x_{1},x_{2})-(x_{2}z)^{-1}\delta\left(\frac{x_{2}-x_{1}}{-x_{2}z}\right)
B(x_{2},x_{1})\nonumber\\
&&\hspace{2cm}=x_{1}^{-1}\delta\left(\frac{x_{2}(1+z)}{x_{1}}\right)
C(\log(1+z),x_{2}).
\end{eqnarray}
Furthermore, the converse is also true. \el

\begin{proof} We have the standard delta-function identity
\begin{eqnarray*}\label{edelta-identity-old}
&&x_{0}^{-1}\delta\left(\frac{x_{1}-x_{2}}{x_{0}}\right)
-x_{0}^{-1}\delta\left(\frac{x_{2}-x_{1}}{-x_{0}}\right)
=x_{1}^{-1}\delta\left(\frac{x_{2}+x_{0}}{x_{1}}\right)
\end{eqnarray*}
(see \cite{flm}). Substituting $x_{0}=x_{2}z$ with $z$ a new formal
variable, we have
\begin{eqnarray}\label{edelta-identity-0}
&&(x_{2}z)^{-1}\delta\left(\frac{x_{1}-x_{2}}{x_{2}z}\right)
-(x_{2}z)^{-1}\delta\left(\frac{x_{2}-x_{1}}{-x_{2}z}\right)
=x_{1}^{-1}\delta\left(\frac{x_{2}(1+z)}{x_{1}}\right),
\end{eqnarray}
where it is understood that  for $n\in \Z$,
$$(1+z)^{n}=\sum_{j\ge 0}\binom{n}{j}z^{j}\in \C[[z]].$$
Then using Remark \ref{rlog-exp} we obtain
\begin{eqnarray*}
&&(x_{2}z)^{-1}\delta\left(\frac{x_{1}-x_{2}}{x_{2}z}\right)(x_{2}z)^{k}
A(x_{1},x_{2})-(x_{2}z)^{-1}\delta\left(\frac{x_{2}-x_{1}}{-x_{2}z}\right)
(x_{2}z)^{k}B(x_{2},x_{1})\nonumber\\
&=&(x_{2}z)^{-1}\delta\left(\frac{x_{1}-x_{2}}{x_{2}z}\right)(x_{1}-x_{2})^{k}
A(x_{1},x_{2})\\
&&\hspace{2cm}-(x_{2}z)^{-1}\delta\left(\frac{x_{2}-x_{1}}{-x_{2}z}\right)
(x_{1}-x_{2})^{k}B(x_{2},x_{1})\nonumber\\
&=&x_{1}^{-1}\delta\left(\frac{x_{2}(1+z)}{x_{1}}\right)
\left((x_{1}-x_{2})^{k} A(x_{1},x_{2})\right)\nonumber\\
&=&x_{1}^{-1}\delta\left(\frac{x_{2}(1+z)}{x_{1}}\right)
\left((x_{1}-x_{2})^{k} A(x_{1},x_{2})\right)|_{x_{1}=x_{2}(1+z)}\nonumber\\
&=&x_{1}^{-1}\delta\left(\frac{x_{2}(1+z)}{x_{1}}\right)
\left(\left((x_{1}-x_{2})^{k} A(x_{1},x_{2})\right)|_{x_{1}=x_{2}e^{x_{0}}}
\right)|_{x_{0}=\log(1+z)}\nonumber\\
&=&x_{1}^{-1}\delta\left(\frac{x_{2}(1+z)}{x_{1}}\right)
(x_{2}e^{x_{0}}-1)^{k}C(x_{0},x_{2})|_{x_{0}=\log(1+z)}\nonumber\\
&=&x_{1}^{-1}\delta\left(\frac{x_{2}(1+z)}{x_{1}}\right)
(x_{2}z)^{k}C(\log(1+z),x_{2}),
\end{eqnarray*}
which implies (\ref{edelta-identity}).

Conversely, assume (\ref{edelta-identity}). Let $k\in \N$ be such
that
$$x_{0}^{k}C(x_{0},x_{2})\in \left(\Hom
(W,W((x_{2})))\right)[[x_{0}]].$$ Then
$$z^{k}C(\log(1+z),x_{2})\in \left(\Hom
(W,W((x_{2})))\right)[[z]],$$ as $\log(1+z)=zg(z)$ with $g(z)\in
\C[[z]]$ invertible. Applying $\Res_{z}z^{k}$ to
(\ref{edelta-identity}) we get
$$(x_{1}-x_{2})^{k}A(x_{1},x_{2})=(x_{1}-x_{2})^{k}B(x_{1},x_{2}).$$
Using this and (\ref{edelta-identity-0}) we get
\begin{eqnarray*}
&&x_{1}^{-1}\delta\left(\frac{x_{2}(1+z)}{x_{1}}\right)
\left((x_{1}-x_{2})^{k}
A(x_{1},x_{2})\right)\\
&=&x_{1}^{-1}\delta\left(\frac{x_{2}(1+z)}{x_{1}}\right)
(x_{2}z)^{k}C(\log(1+z),x_{2}).
\end{eqnarray*}
Substituting $z=e^{x_{0}}-1$, we get
\begin{eqnarray*}
x_{1}^{-1}\delta\left(\frac{x_{2}e^{x_{0}}}{x_{1}}\right)
\left((x_{1}-x_{2})^{k} A(x_{1},x_{2})\right)
=x_{1}^{-1}\delta\left(\frac{x_{2}e^{x_{0}}}{x_{1}}\right)
(x_{2}z)^{k}C(x_{0},x_{2}),
\end{eqnarray*}
which implies
$$\left((x_{1}-x_{2})^{k}A(x_{1},x_{2})\right)|_{x_{1}=x_{2}e^{x_{0}}}
=x_{2}^{k}(e^{x_{0}}-1)^{k}C(x_{0},x_{2}).$$ This completes the
proof.
\end{proof}

Now we are ready to present our second main result of this section.

\bp{pwqva-phimodule} Let $V$ be a weak quantum vertex algebra and
let $(W,Y_{W})$ be a $\phi$-coordinated module for $V$ viewed as a
nonlocal vertex algebra. Let $u,v\in V$ and assume that
$$(x_{1}-x_{2})^{k}Y(u,x_{1})Y(v,x_{2})=(x_{1}-x_{2})^{k}
\sum_{i=1}^{r}\iota_{x_{2},x_{1}}(f_{i}(e^{x_{1}-x_{2}}))
Y(v^{(i)},x_{2})Y(u^{(i)},x_{1})$$ with $k\in \N,\; f_{i}(x)\in
\C(x),\; u^{(i)},v^{(i)}\in V$ for $1\le i\le r$. Then
\begin{eqnarray}\label{ejacobi-new4}
&&(x_{2}z)^{-1}\delta\left(\frac{x_{1}-x_{2}}{x_{2}z}\right)
Y_{W}(u,x_{1})Y_{W}(v,x_{2})\nonumber\\
&&\ \ \ \
-(x_{2}z)^{-1}\delta\left(\frac{x_{2}-x_{1}}{-x_{2}z}\right)
\sum_{i=1}^{r}\iota_{x_{2},x_{1}}(f_{i}(x_{1}/x_{2}))
Y_{W}(v^{(i)},x_{2})Y_{W}(u^{(i)},x_{1})\nonumber\\
&=&x_{1}^{-1}\delta\left(\frac{x_{2}(1+z)}{x_{1}}\right)
Y_{W}(Y(u,\log(1+z))v,x_{2}).
\end{eqnarray}
Furthermore, we have
\begin{eqnarray}\label{ejacobi-new5}
&&Y_{W}(u,x_{1})Y_{W}(v,x_{2})
-\sum_{i=1}^{r}\iota_{x_{2},x_{1}}(f_{i}(x_{1}/x_{2}))
Y_{W}(v^{(i)},x_{2})Y_{W}(u^{(i)},x_{1})\nonumber\\
&=&\Res_{x_{0}}x_{1}^{-1}\delta\left(\frac{x_{2}e^{x_{0}}}{x_{1}}\right)
x_{2}e^{x_{0}}Y_{W}(Y(u,x_{0})v,x_{2}).
\end{eqnarray}
 \ep

\begin{proof} From definition,
there exists a nonnegative integer $l$ such that
$$(x_{1}-x_{2})^{l}Y_{W}(u,x_{1})Y_{W}(v,x_{2})\in \Hom (W,W((x_{1},x_{2})))$$
and
$$x_{2}^{l}(e^{x_{0}}-1)^{l}Y_{W}(Y(u,x_{0})v,x_{2})
=\left((x_{1}-x_{2})^{l}Y_{W}(u,x_{1})Y_{W}(v,x_{2})\right)|_{x_{1}=x_{2}e^{x_{0}}}.$$
On the other hand, by Proposition \ref{palgebra-module} we also have
\begin{eqnarray*}
&&(x_{1}-x_{2})^{l}Y_{W}(u,x_{1})Y_{W}(v,x_{2})\\
&=&(x_{1}-x_{2})^{l}
\sum_{i=1}^{r}\iota_{x_{2},x_{1}}(f_{i}(x_{1}/x_{2}))
Y_{W}(v^{(i)},x_{2})Y_{W}(u^{(i)},x_{1}).
\end{eqnarray*}
Then the first assertion follows immediately from Lemma
\ref{ldelta-identity}. Furthermore, applying $\Res_{z}x_{2}$ we get
\begin{eqnarray*}
&&Y_{W}(u,x_{1})Y_{W}(v,x_{2})
-\sum_{i=1}^{r}\iota_{x_{2},x_{1}}(f_{i}(x_{1}/x_{2}))
Y_{W}(v^{(i)},x_{2})Y_{W}(u^{(i)},x_{1})\nonumber\\
&=&\Res_{z}x_{1}^{-1}\delta\left(\frac{x_{2}(1+z)}{x_{1}}\right)
x_{2}Y_{W}(Y(u,\log(1+z))v,x_{2})\\
&=&\Res_{x_{0}}x_{1}^{-1}\delta\left(\frac{x_{2}e^{x_{0}}}{x_{1}}\right)
x_{2}e^{x_{0}}Y_{W}(Y(u,x_{0})v,x_{2}),
\end{eqnarray*}
proving the second assertion. \end{proof}

As a consequence we have:

\bc{cdiff-expression} Let $W$ be a vector space and let $U$ be an
$\S_{trig}$-local subset of $\E(W)$. Assume that
$$a(x),b(x),a^{(i)}(x),b^{(i)}(x)\in U,\; f_{i}(x)\in \C(x) \ \
(i=1,\dots,r)$$ satisfy
$$(x_{1}-x_{2})^{k}a(x_{1})b(x_{2})
=(x_{1}-x_{2})^{k}\sum_{i=1}^{r}\iota_{x_{2},x_{1}}(f_{i}(x_{1}/x_{2}))b^{(i)}(x_{2})a^{(i)}(x_{1})$$
for some nonnegative integer $k$. Then
\begin{eqnarray}\label{ejacobi-new3}
&&(xz)^{-1}\delta\left(\frac{x_{1}-x}{xz}\right)
a(x_{1})b(x)-(xz)^{-1}\delta\left(\frac{x-x_{1}}{-xz}\right)
\sum_{i=1}^{r}f_{i}(x_{1}/x)b^{(i)}(x)a^{(i)}(x_{1})\nonumber\\
&=&x_{1}^{-1}\delta\left(\frac{x(1+z)}{x_{1}}\right)
Y_{\E}^{e}(a(x),\log(1+z))b(x).
\end{eqnarray}
In particular,
\begin{eqnarray}\label{edef-future}
&&Y_{\E}^{e}(a(x),\log(1+z))b(x)
=\Res_{x_{1}}(xz)^{-1}\delta\left(\frac{x_{1}-x}{xz}\right)
a(x_{1})b(x)\nonumber\\
&&\hspace{2cm}-\Res_{x_{1}}(xz)^{-1}\delta\left(\frac{x-x_{1}}{-xz}\right)
\sum_{i=1}^{r}f_{i}(x_{1}/x)b^{(i)}(x)a^{(i)}(x_{1}).
\end{eqnarray} \ec

\begin{proof} As $U$ is $\S_{trig}$-local, by Theorem \ref{tqva}
$\<U\>_{e}$ is a weak quantum vertex algebra with $W$ as a
$\phi$-coordinated module. By Proposition \ref{pwqva-phimodule}, we
have (\ref{ejacobi-new3}), which immediately implies
(\ref{edef-future}).
\end{proof}

\section{Quantum vertex algebras associated with quantum $\beta\gamma$-system}
In this section we associate quantum vertex algebras to certain
quantum $\beta\gamma$-systems, using the general machinery developed
in previous sections. We first introduce a quantum
$\beta\gamma$-system of trigonometric type, which is a modified
version of the quantum $\beta\gamma$-system in \cite{efk}. Then we
introduce a quantum $\beta\gamma$-system of rational type, to
describe the quantum vertex algebras constructed from the
trigonometric type quantum $\beta\gamma$-system.

We start with the following quantum $\beta\gamma$-system:

\bd{dqbg-trig} {\em Let $q$ be a nonzero complex number. Define
$A_{q}^{trig}(\beta\gamma)$ to be the associative algebra with
identity over $\C$ with generators $\tilde{\beta}_{n},\;
\tilde{\gamma}_{n}\; (n\in \Z)$, which are considered as the
coefficients of the currents
$$ \tilde{\beta}(x)=\sum_{n\in \Z}\tilde{\beta}_{n}x^{-n-1},\ \ \ \
  \tilde{\gamma}(x)=\sum_{n\in \Z}\tilde{\gamma}_{n}x^{-n-1},$$
subject to relations
\begin{eqnarray}\label{etrig-qbg-1}
&&\tilde{\beta}(x)\tilde{\beta}(z)
=\iota_{z,x}\left(\frac{x-qz}{qx-z}\right)
\tilde{\beta}(z)\tilde{\beta}(x),
\nonumber\\
&&\tilde{\gamma}(x)\tilde{\gamma}(z)
=\iota_{z,x}\left(\frac{x-qz}{qx-z}\right)\tilde{\gamma}(z)\tilde{\gamma}(x),
\nonumber\\
&&\tilde{\beta}(x)\tilde{\gamma}(z)
-\iota_{z,x}\left(\frac{qx-z}{x-qz}\right)\tilde{\gamma}(z)\tilde{\beta}(x)
=\delta\left(\frac{x}{z}\right).
\end{eqnarray}}
\ed

When $q=1$, it can be readily seen that  $\tilde{\beta}(x)$ and
$x^{-1}\tilde{\gamma}(x)$ form the standard $\beta\gamma$-system
where $A_{q}^{trig}(\beta\gamma)$ is isomorphic to the universal
enveloping algebra of an infinite-dimensional Heisenberg Lie
algebra.

\br{rdefinition}{\em Here we give some details about the definition
of $A_{q}^{trig}(\beta\gamma)$.  Set
$$\iota_{z,x}\left(\frac{x-qz}{qx-z}\right)= \sum_{k\ge 0}\lambda_{k}(x/z)^{k},\ \
\iota_{z,x}\left(\frac{qx-z}{x-qz}\right)=\sum_{k\ge
0}\lambda'_{k}(x/z)^{k}$$ with $\lambda_{k},\lambda_{k}'\in \C$ for
$k\ge 0$. The defining relations (\ref{etrig-qbg-1}) amount to
\begin{eqnarray}\label{etrig-qbg-comp}
&&\tilde{\beta}_{m}\tilde{\beta}_{n} =\sum_{k\ge
0}\lambda_{k}\tilde{\beta}_{n-k}\tilde{\beta}_{m+k},\ \ \ \
\tilde{\gamma}_{m}\tilde{\gamma}_{n} =\sum_{k\ge
0}\lambda_{k}\tilde{\gamma}_{n-k}\tilde{\gamma}_{m+k},\nonumber\\
&&\hspace{1cm}\tilde{\beta}_{m}\tilde{\gamma}_{n} -\sum_{k\ge
0}\lambda'_{k}\tilde{\gamma}_{n-k}\tilde{\beta}_{m+k}=\delta_{m+n+2,0}
\end{eqnarray}
for $m,n\in \Z$. Let $T$ be the free associative algebra over $\C$
with generators $\beta_{n},\; \gamma_{n}$ ($n\in \Z$). Define
$$\deg \beta_{n}=\deg \gamma_{n}=n+1
\ \ \ \mbox{ for }n\in \Z,$$ to make $T$ a $\Z$-graded algebra whose
homogeneous subspace of degree-$n$ is denoted by $T(n)$. Following
\cite{fz}, for $n\in \Z,\; k\ge 1$, set
$$T(n,k)=\sum_{r\ge k} T(n-r)T(r)\subset T(n).$$
We have $\cap_{k\ge 1}T(n,k)=0$. Equip $T(n)$ with the topology with
$a+T(n,k)$ for $a\in T(n),\; k\ge 0$ as a basis of open sets. Let
$\overline{T(n)}$ be the formal completion of $T(n)$. Set
$\overline{T}=\oplus_{n\in \Z}\overline{T(n)}$. Then the algebra
$A_{q}^{trig}(\beta\gamma)$ can be defined as the quotient algebra
of $\overline{T}$ modulo the relations (\ref{etrig-qbg-comp}). Since
all the relations are homogeneous, $A_{q}^{trig}(\beta\gamma)$ is a
$\Z$-graded algebra.} \er

As $(x-z)\delta\left(\frac{x}{z}\right)=0$, from the third relation
in (\ref{etrig-qbg-1}) we get
\begin{eqnarray}\label{ebeta-gamma-trig}
(x-z)\tilde{\beta}(x)\tilde{\gamma}(z)
=(x-z)\iota_{z,x}\left(\frac{qx-z}{x-qz}\right)\tilde{\gamma}(z)\tilde{\beta}(x).
\end{eqnarray}

By a {\em restricted} $A_{q}^{trig}(\beta\gamma)$-module we mean an
$A_{q}^{trig}(\beta\gamma)$-module $W$ which equipped with the
discrete topology is a continuous module. Then a restricted
$A_{q}^{trig}(\beta\gamma)$-module amounts to a $T$-module $W$ such
that for every $w\in W$, $\beta_{n}w=\gamma_{n}w=0$ for $n$
sufficiently large and the relations corresponding to
(\ref{etrig-qbg-comp}) applied to $w$ hold. Let $W$ be a restricted
$A_{q}^{trig}(\beta\gamma)$-module. With the relations
(\ref{etrig-qbg-1}) and (\ref{ebeta-gamma-trig}), we see that $\{
\tilde{\beta}(x),\tilde{\gamma}(x)\}$ is an $\S_{trig}$-local subset
of $\E(W)$. In view of Theorem \ref{tqva}, $\{
\tilde{\beta}(x),\tilde{\gamma}(x)\}$ generates a weak quantum
vertex algebra $V_{W}$ inside $\E(W)$. To describe the structure of
$V_{W}$ we need another algebra.

\bd{dqbg-rational-1} {\em Let $q$ be a complex number. Define
$A_{q}^{rat}(\beta\gamma)$ to be the associative algebra with
identity over $\C$ with generators $\hat{\beta}_{n},\;
\hat{\gamma}_{n}\; (n\in \Z)$, subject to relations
\begin{eqnarray}\label{erat-qbg-1}
&&\hat{\beta}(x)\hat{\beta}(z)
=\iota_{z,x}\left(\frac{e^{x-z}-q}{qe^{x-z}-1}\right)\hat{\beta}(z)\hat{\beta}(x),
\nonumber\\
&&\hat{\gamma}(x)\hat{\gamma}(z)
=\iota_{z,x}\left(\frac{e^{x-z}-q}{qe^{x-z}-1}\right)\hat{\gamma}(z)\hat{\gamma}(x),
\nonumber\\
&&\hat{\beta}(x)\hat{\gamma}(z)
-\iota_{z,x}\left(\frac{qe^{x-z}-1}{e^{x-z}-q}\right)\hat{\gamma}(z)\hat{\beta}(x)
=z^{-1}\delta\left(\frac{x}{z}\right).
\end{eqnarray}}
\ed

\br{rdef-rational} {\em  Here is a precise definition of the algebra
$A_{q}^{rat}(\beta\gamma)$. Set
$$\frac{e^{x-z}-q}{qe^{x-z}-1}=\sum_{k\ge 0}\mu_{k}(x-z)^{k},\ \
\ \ \frac{qe^{x-z}-1}{e^{x-z}-q}=\sum_{k\ge 0}\mu'_{k}(x-z)^{k}$$
with $\mu_{k},\mu_{k}'\in \C$ for $k\ge 0$. The defining relations
(\ref{erat-qbg-1}) read as
\begin{eqnarray}\label{erat-qbg-comp}
&&\hat{\beta}_{m}\hat{\beta}_{n} =\sum_{k,i\ge
0}\binom{k}{i}(-1)^{i}\mu_{k}\hat{\beta}_{n+i}\hat{\beta}_{m+k-i}, \
\ \ \hat{\gamma}_{m}\hat{\gamma}_{n} =\sum_{k,i\ge
0}\binom{k}{i}(-1)^{i}\mu_{k}\hat{\gamma}_{n+i}\hat{\gamma}_{m+k-i},
\nonumber\\
&&\hspace{1cm}\hat{\beta}_{m}\hat{\gamma}_{n} -\sum_{k,i\ge
0}\binom{k}{i}(-1)^{i}\mu'_{k}\hat{\gamma}_{n+i}\hat{\beta}_{m+k-i}=\delta_{m+n+1,0}
\end{eqnarray}
for $m,n\in \Z$.  Let $T$ be the free associative algebra  as in
Remark \ref{rdefinition}, generated by $\beta_{n},\gamma_{n}$ for
$n\in \Z$, and equip $T$ with the same $\Z$-grading. For $k\ge 0$,
set
$$T[k]=\sum_{n\ge k}T(n)\subset T.$$
Equip $T$ with the topology with $a+T(k)$ for $a\in T,\; k\ge 0$ as
a basis of open sets. Let $\tilde{T}$ be the formal completion of
$T$. The algebra $A_{q}^{rat}(\beta\gamma)$ can be defined as the
quotient algebra of $\tilde{T}$ modulo the relations
(\ref{erat-qbg-comp}). Since the defining relations are
inhomogeneous, the algebra $A_{q}^{rat}(\beta\gamma)$ is not
 $\Z$-graded in the obvious way. } \er

\br{rspecial-cases} {\em Notice that when $q=1$, the quantum
$\beta\gamma$-system defined in Definition \ref{dqbg-rational-1} is
exactly the standard $\beta\gamma$-system. If $q=-1$, the defining
relations become
\begin{eqnarray}\label{enew-rat-q=-1}
&&\hat{\beta}(x)\hat{\beta}(z) =-\hat{\beta}(z)\hat{\beta}(x),
\nonumber\\
&&\hat{\gamma}(x)\hat{\gamma}(z) =-\hat{\gamma}(z)\hat{\gamma}(x),
\nonumber\\
&&\hat{\beta}(x)\hat{\gamma}(z) +\hat{\gamma}(z)\hat{\beta}(x)
=z^{-1}\delta\left(\frac{x}{z}\right).
\end{eqnarray}
In this case, $A_{q}^{rat}(\beta\gamma)$ is an
(infinite-dimensional) Clifford algebra.} \er

Just as with $A_{q}^{trig}(\beta\gamma)$, by a {\em restricted
$A_{q}^{rat}(\beta\gamma)$-module} we mean an
$A_{q}^{rat}(\beta\gamma)$-module $W$ which equipped with the
discrete topology is a continuous module. A restricted
$A_{q}^{rat}(\beta\gamma)$-module simply amounts to a module for the
free algebra $T$ such that for any $w\in W$,
$\beta_{n}w=0=\gamma_{n}w$ for $n$ sufficiently large and such that
the relations corresponding to (\ref{erat-qbg-comp}) after applied
to each vector $w\in W$ hold.

\bd{dvacuum-module} {\em A {\em vacuum
$A_{q}^{rat}(\beta\gamma)$-module} is a restricted
 $A_{q}^{rat}(\beta\gamma)$-module $W$
equipped with a vector $w_{0}\in W$, satisfying the condition that
 $W=A_{q}^{rat}(\beta\gamma)w_{0}$,
 \begin{eqnarray}
\hat{\beta}_{n}w_{0}=\hat{\gamma}_{n}w_{0}=0 \ \ \ \ \mbox{ for
}n\ge 0.
\end{eqnarray}
We sometimes denote a vacuum module by a pair $(W,w_{0})$. } \ed

We are going to prove that the weak quantum vertex algebra $V_{W}$
associated to a restricted $A_{q}^{trig}(\beta\gamma)$-module $W$ is
naturally a vacuum $A_{q}^{rat}(\beta\gamma)$-module. To achieve
this goal, we shall need the following technical result:

\bl{labuv} Let $W$ be a vector space and let $a(x),b(x)\in \E(W)$.
Assume that there exist
$$0\ne p(x)\in
\C[x],\; q_{i}(x)\in \C((x)),\;u^{(i)}(x),v^{(i)}(x)\in \E(W) \ \
(1\le i\le r)$$ such that
\begin{eqnarray}\label{epab=quv0}
p(x_{1}/x_{2})a(x_{1})b(x_{2})=\sum_{i=1}^{r}q_{i}(x_{1}/x_{2})u^{(i)}(x_{2})v^{(i)}(x_{1}).
\end{eqnarray}
Then $(a(x),b(x))$ is quasi compatible and
\begin{eqnarray*}
&&p(e^{x_{0}})Y_{\E}^{e}(a(x),x_{0})b(x)\\
&=&\Res_{x_{1}}\left(\frac{1}{x_{1}-xe^{x_{0}}}p(x_{1}/x)a(x_{1})b(x)-
\frac{1}{-xe^{x_{0}}+x_{1}}\sum_{i=1}^{r}q_{i}(x_{1}/x)u^{(i)}(x)v^{(i)}(x_{1})\right).
\end{eqnarray*}
Furthermore, if $k$ is the order of zero of $p(x)$ at $1$, then
$a(x)_{n}^{e}b(x)=0$ for $n\ge k$ and
\begin{eqnarray}\label{ea-1b}
&&\frac{1}{k!}p^{(k)}(1)a(x)_{k-1}^{e}b(x)\nonumber\\
&=&\Res_{x_{1}}\left(\frac{1}{x_{1}-x}p(x_{1}/x)a(x_{1})b(x)-
\frac{1}{-x+x_{1}}q(x_{1}/x)u^{(i)}(x)v^{(i)}(x_{1})\right).\ \
\end{eqnarray}
\el

\begin{proof} By observing both sides of (\ref{epab=quv0}) we see
that
$$p(x_{1}/x_{2})a(x_{1})b(x_{2})\in \Hom (W,W((x_{1},x_{2}))),$$
which implies that $(a(x),b(x))$ is quasi compatible. Furthermore,
we have
\begin{eqnarray*}
&&p(e^{x_{0}})Y_{\E}^{e}(a(x),x_{0})b(x)\\
&=&\left(p(x_{1}/x)a(x_{1})b(x)\right)|_{x_{1}=xe^{x_{0}}}\\
&=&\Res_{x_{1}}x_{1}^{-1}\delta\left(\frac{xe^{x_{0}}}{x_{1}}\right)
\left(p(x_{1}/x)a(x_{1})b(x)\right)\\
&=&\Res_{x_{1}}\left(\frac{1}{x_{1}-xe^{x_{0}}}p(x_{1}/x)a(x_{1})b(x)-
\frac{1}{-xe^{x_{0}}+x_{1}}\left(p(x_{1}/x)a(x_{1})b(x)\right)\right)\\
&=&\Res_{x_{1}}\left(\frac{1}{x_{1}-xe^{x_{0}}}p(x_{1}/x)a(x_{1})b(x)-
\frac{1}{-xe^{x_{0}}+x_{1}}\sum_{i=1}^{r}q_{i}(x_{1}/x)u^{(i)}(x)v^{(i)}(x_{1})\right)
\end{eqnarray*}
as $x_{1}^{-1}\delta\left(\frac{xe^{x_{0}}}{x_{1}}\right)
=\frac{1}{x_{1}-xe^{x_{0}}}-\frac{1}{-xe^{x_{0}}+x_{1}}$. This
proves the first part of the lemma.

Note that $p(e^{x_{0}})Y_{\E}^{e}(a(x),x_{0})b(x)$ involves only
nonnegative powers of $x_{0}$. As $k$ is the order of zero of $p(x)$
at $1$, we have $p(e^{x_{0}})=x_{0}^{k}g(x_{0})$ for some $g(x)\in
\C[[x]]$ with $g(0)\ne 0$. Since $g(x_{0})$ is a unit in
$\C[[x_{0}]]$, we have that $x_{0}^{k}Y_{\E}^{e}(a(x),x_{0})b(x)$
involves only nonnegative powers of $x_{0}$. That is,
$a(x)_{n}^{e}b(x)=0$ for $n\ge k$. Then, applying
$\Res_{x_{0}}x_{0}^{-1}$, (or setting $x_{0}=0$), we obtain
(\ref{ea-1b}).
\end{proof}

Now we have:

\bp{pbeta-gamma-vacuum} Let $W$ be a restricted
$A_{q}^{trig}(\beta\gamma)$-module and let $V_{W}$ be the weak
quantum vertex algebra generated by the $\S_{trig}$-local subset
$\{\tilde{\beta}(x),\tilde{\gamma}(x)\}$ of $\E(W)$. Then $V_{W}$ is
an $A_{q}^{rat}(\beta\gamma)$-module with $\hat{\beta}(z)$ and
$\hat{\gamma}(z)$ acting as $Y_{\E}^{e}(\tilde{\beta}(x),z)$ and
$Y_{\E}^{e}(\tilde{\gamma}(x),z)$, respectively. Furthermore,
$(V_{W},1_{W})$ is a vacuum $A_{q}^{rat}(\beta\gamma)$-module. \ep

\begin{proof} With the relations (\ref{etrig-qbg-1}) and (\ref{ebeta-gamma-trig}),
in view of Proposition \ref{pconvert} we have
\begin{eqnarray*}
&&Y_{\E}^{e}(\tilde{\beta}(x),x_{1})Y_{\E}^{e}(\tilde{\beta}(x),x_{2})
=\left(\frac{e^{x_{1}-x_{2}}-q}{qe^{x_{1}-x_{2}}-1}\right)
Y_{\E}^{e}(\tilde{\beta}(x),x_{2})Y_{\E}^{e}(\tilde{\beta}(x),x_{1}),\\
&&Y_{\E}^{e}(\tilde{\gamma}(x),x_{1})Y_{\E}^{e}(\tilde{\gamma}(x),x_{2})
=\left(\frac{e^{x_{1}-x_{2}}-q}{qe^{x_{1}-x_{2}}-1}\right)
Y_{\E}^{e}(\tilde{\gamma}(x),x_{2})Y_{\E}^{e}(\tilde{\gamma}(x),x_{1}),\\
&&(x_{1}-x_{2})Y_{\E}^{e}(\tilde{\beta}(x),x_{1})Y_{\E}^{e}(\tilde{\gamma}(x),x_{2})\\
&&\hspace{2cm}=(x_{1}-x_{2})\left(\frac{qe^{x_{1}-x_{2}}-1}{e^{x_{1}-x_{2}}-q}\right)
Y_{\E}^{e}(\tilde{\gamma}(x),x_{2})Y_{\E}^{e}(\tilde{\beta}(x),x_{1}).
\end{eqnarray*}
Furthermore, due to the last relation, we have
\begin{eqnarray}\label{eproof-sjacobi-proof}
&&x_{0}^{-1}\delta\left(\frac{x_{1}-x_{2}}{x_{0}}\right)
Y_{\E}^{e}(\tilde{\beta}(x),x_{1})Y_{\E}^{e}(\tilde{\gamma}(x),x_{2})\nonumber\\
&&\hspace{1.5cm}
-x_{0}^{-1}\delta\left(\frac{x_{2}-x_{1}}{-x_{0}}\right)
\left(\frac{qe^{x_{1}-x_{2}}-1}{e^{x_{1}-x_{2}}-q}\right)
Y_{\E}^{e}(\tilde{\gamma}(x),x_{2})Y_{\E}^{e}(\tilde{\beta}(x),x_{1})\nonumber\\
&=&x_{2}^{-1}\delta\left(\frac{x_{1}-x_{0}}{x_{2}}\right)
Y_{\E}^{e}(Y_{\E}^{e}(\tilde{\beta}(x),x_{0})\tilde{\gamma}(x),x_{2}).
\end{eqnarray}
Combining (\ref{ebeta-gamma-trig}) with Lemma \ref{labuv} we get
$\tilde{\beta}(x)_{n}^{e}\tilde{\gamma}(x)=0$ for $n\ge 1$ and we
have
\begin{eqnarray*}
\tilde{\beta}(x)_{0}^{e}\tilde{\gamma}(x)
&=&\Res_{x_{1}}\left(x^{-1}\tilde{\beta}(x_{1})\tilde{\gamma}(x)-x^{-1}
\left(\frac{qx_{1}/x-1}{x_{1}/x-q}\right)\tilde{\gamma}(x)\tilde{\beta}(x_{1})\right)\\
&=&\Res_{x_{1}}x^{-1}\delta\left(\frac{x_{1}}{x}\right)\\
&=&1.
\end{eqnarray*}
Then applying $\Res_{x_{0}}$ to (\ref{eproof-sjacobi-proof}) we
obtain
\begin{eqnarray*}
&&Y_{\E}^{e}(\tilde{\beta}(x),x_{1})Y_{\E}^{e}(\tilde{\gamma}(x),x_{2})
- \left(\frac{qe^{x_{1}-x_{2}}-1}{e^{x_{1}-x_{2}}-q}\right)
Y_{\E}^{e}(\tilde{\gamma}(x),x_{2})Y_{\E}^{e}(\tilde{\beta}(x),x_{1})\\
&=&x_{2}^{-1}\delta\left(\frac{x_{1}}{x_{2}}\right)
Y_{\E}^{e}(\tilde{\beta}(x)_{0}^{e}\tilde{\gamma}(x),x_{2})\\
&=&x_{2}^{-1}\delta\left(\frac{x_{1}}{x_{2}}\right).
\end{eqnarray*}
Now, we see that with $\hat{\beta}(z)$ and $\hat{\gamma}(z)$ acting
as $Y_{\E}^{e}(\tilde{\beta}(x),z)$ and
$Y_{\E}^{e}(\tilde{\gamma}(x),z)$, respectively, $V_{W}$ becomes an
$A_{q}^{rat}(\beta\gamma)$-module. Since  $V_{W}$ as a nonlocal
vertex algebra is generated by $\hat{\beta}(x)$ and
$\hat{\gamma}(x)$, it follows that $V_{W}$ as an
$A_{q}^{rat}(\beta\gamma)$-module is generated by $1_{W}$. We have
$\tilde{\beta}(x)_{n}1_{W}=0$ and $\tilde{\gamma}(x)_{n}1_{W}=0$ for
$n\ge 0$. Therefore $(V_{W},1_{W})$ is a vacuum
$A_{q}^{rat}(\beta\gamma)$-module.
\end{proof}

Next, we construct a universal vacuum
$A_{q}^{rat}(\beta\gamma)$-module, following \cite{li-qva2} (Section
4). Let $T_{+}$ denote the subspace of $T$, linearly spanned by the
vectors
$$a^{(1)}_{n_{1}}\cdots a^{(r)}_{n_{r}}$$
for $r\ge 1,\; a^{(i)}\in \{\beta,\gamma\},\; n_{i}\in \Z$ with
$n_{1}+\cdots +n_{r}\ge 0$. As $TT_{+}$ is a left ideal of $T$,
$T/TT_{+}$ is naturally a $T$-module, which we denote by
$\widetilde{V}(\beta\gamma)$. From definition, for any $w\in
\widetilde{V}(\beta\gamma)$, $T(n)w=0$ for $n$ sufficiently large.
Then we let $J_{q}(\beta\gamma)$ be the submodule of
$\widetilde{V}(\beta\gamma)$, generated by the following vectors:
\begin{eqnarray*}
&&\beta_{m}\beta_{n}w-\sum_{k,i\ge
0}\binom{k}{i}(-1)^{i}\mu_{k}\beta_{n+i}\beta_{m+k-i}w, \\
&&\gamma_{m}\gamma_{n}w-\sum_{k,i\ge
0}\binom{k}{i}(-1)^{i}\mu_{k}\gamma_{n+i}\gamma_{m+k-i}w,
\nonumber\\
&&\beta_{m}\gamma_{n}w -\sum_{k,i\ge
0}\binom{k}{i}(-1)^{i}\mu'_{k}\gamma_{n+i}\beta_{m+k-i}w-\delta_{m+n+1,0}w
\end{eqnarray*}
for $m,n\in \Z,\; w\in \widetilde{V}(\beta\gamma)$ (recall Remark
\ref{rdef-rational}). Set
\begin{eqnarray}
V_{q}(\beta\gamma)=\widetilde{V}(\beta\gamma)/J_{q}(\beta\gamma)
\end{eqnarray}
and set ${\bf 1}=1+J_{q}(\beta\gamma)\in V_{q}(\beta\gamma)$. {}From
the construction, $(V_{q}(\beta\gamma),{\bf 1})$ is naturally a
vacuum $A_{q}^{rat}(\beta\gamma)$-module. Set
\begin{eqnarray}
\hat{\beta}=\hat{\beta}_{-1}{\bf 1},\ \
\hat{\gamma}=\hat{\gamma}_{-1}{\bf 1} \in V_{q}(\beta\gamma).
\end{eqnarray}

\bt{tbeta-gamma-qva} Let $q$ be any nonzero complex number. The
vacuum $A_{q}^{rat}(\beta\gamma)$-module $(V_{q}(\beta\gamma),{\bf
1})$ is universal in the obvious sense  and there exists a weak
quantum vertex algebra structure on $V_{q}(\beta\gamma)$, which is
uniquely determined by the condition that ${\bf 1}$ is the vacuum
vector and
$$Y(\hat{\beta},x)=\hat{\beta}(x),\;\;
Y(\hat{\gamma},x)=\hat{\gamma}(x).$$ Furthermore,
$V_{q}(\beta\gamma)$ is an irreducible quantum vertex algebra. On
the other hand, for every restricted
$A_{q}^{rat}(\beta\gamma)$-module $W$, there exists a
$V_{q}(\beta\gamma)$-module structure $Y_{W}$ on $W$, which is
uniquely determined by the condition that
$$Y_{W}(\hat{\beta},x)=\hat{\beta}(x),\;\;
Y_{W}(\hat{\gamma},x)=\hat{\gamma}(x).$$ \et

\begin{proof}  Let $H$ be a
vector space with $\{\beta,\gamma\}$ as a basis and define a linear
map $\S(x): H\otimes H\rightarrow H\otimes H\otimes \C[[x]]$ by
\begin{eqnarray*}
&&\S(x)(\beta \otimes \beta)=(\beta \otimes \beta)f(x),\ \ \
\S(x)(\gamma \otimes \gamma)=(\gamma \otimes
\gamma)f(x),\\
&&\S(x)(\beta \otimes \gamma)=(\beta \otimes \gamma)f(x),\ \ \
\S(x)(\gamma \otimes \beta)=(\gamma \otimes \beta)g(x),
\end{eqnarray*}
where $f(x)$ and $g(x)$ are the formal Taylor series expansions at
$0$ of $\frac{e^{-x}-q}{qe^{-x}-1}$ and
$\frac{qe^{-x}-1}{e^{-x}-q}$, respectively. Then
$V_{q}(\beta\gamma)$ is simply the $(H,\S)$-module $V(H,\S)$ in
\cite{li-qva2}. In view of this, universality follows from
Proposition 4.3 of \cite{li-qva2}. The assertion on weak quantum
vertex algebra structure and the assertion on module structure
follow immediately from Proposition 4.2 of \cite{kl} (cf.
\cite{li-qva2}, Propositions 2.18 and 4.3). As for the
irreducibility assertion we shall use a result of \cite{kl}. A
special case of Theorem 4.9 of \cite{kl} states that for any
$p(x)\in \C[[x]]$ with $p(0)=1$, there exists a (nonzero) weak
quantum vertex algebra $V$ which is generated by two linearly
independent vectors $u$ and $v$ such that
\begin{eqnarray*}
&&Y(u,x_{1})Y(u,x_{2})=-(p(x_{1}-x_{2})/p(x_{2}-x_{1}))Y(u,x_{2})Y(u,x_{1}),\\
&&Y(v,x_{1})Y(v,x_{2})=-(p(x_{1}-x_{2})/p(x_{2}-x_{1}))Y(v,x_{2})Y(v,x_{1}),\\
&&Y(u,x_{1})Y(v,x_{2})+(p(x_{2}-x_{1})/p(x_{1}-x_{2}))Y(v,x_{2})Y(u,x_{1})
=x_{2}^{-1}\delta\left(\frac{x_{1}}{x_{2}}\right),
\end{eqnarray*}
and furthermore, all such $V$ are irreducible quantum vertex
algebras and isomorphic to each other. For $q=1$, it is known that
$V_{q}(\beta\gamma)$ is a simple (equivalently irreducible) vertex
algebra. Assume $q\ne 1$. Set
$$p(x)=(e^{-x/2}-qe^{x/2})/(1-q)\in \C[[x]].$$ We have $p(0)=1$ and
$$p(-x)/p(x)=\frac{e^{x}-q}{1-qe^{x}}=-\frac{e^{x}-q}{qe^{x}-1}.$$
Then by Theorem 4.9 of \cite{kl}, $V_{q}(\beta\gamma)$ is an
irreducible quantum vertex algebra.
\end{proof}

\br{rvacuum-module} {\em As $V_{q}(\beta\gamma)$ is an irreducible
quantum vertex algebra, $V_{q}(\beta\gamma)$ is an irreducible
$A_{q}^{rat}(\beta\gamma)$-module. It follows that every nonzero
vacuum $A_{q}^{rat}(\beta\gamma)$-module is irreducible and
isomorphic to $V_{q}(\beta\gamma)$. } \er

The following is a connection between quantum vertex algebra
$V_{q}(\beta\gamma)$ and restricted
$A_{q}^{trig}(\beta\gamma)$-modules:

\bt{tbeta-gamma-exp} Let $q$ be a nonzero complex number and let $W$
be a restricted $A_{q}^{trig}(\beta\gamma)$-module. Then there
exists a $\phi$-coordinated $V_{q}(\beta\gamma)$-module structure
$Y_{W}$ on $W$ with $\phi(x,z)=xe^{z}$, which is uniquely determined
by the condition that
$$Y_{W}(\hat{\beta},x)=\tilde{\beta}(x),\;\;
Y_{W}(\hat{\gamma},x)=\tilde{\gamma}(x).$$ On the other hand, for
any $\phi$-coordinated $V_{q}(\beta\gamma)$-module $(W,Y_{W})$, $W$
is a restricted $A_{q}^{trig}(\beta\gamma)$-module with
$\tilde{\beta}(x)$ and $\tilde{\gamma}(x)$ acting as
$Y_{W}(\hat{\beta},x)$ and $Y_{W}(\hat{\gamma},x)$, respectively.
\et

\begin{proof} It is similar to the proof of Theorem \ref{tbeta-gamma-qva}.
First, by Proposition \ref{pbeta-gamma-vacuum}, the weak quantum
vertex algebra $V_{W}$ with $1_{W}$ is a vacuum
$A_{q}^{rat}(\beta\gamma)$-module. As the vacuum
$A_{q}^{rat}(\beta\gamma)$-module $V_{q}(\beta\gamma)$ is universal,
there exists an $A_{q}^{rat}(\beta\gamma)$-module homomorphism
$\theta$ from $V_{q}(\beta\gamma)$ to $V_{W}$, sending ${\bf 1}$ to
$1_{W}$. It follows that $\theta$ is a homomorphism of weak quantum
vertex algebras. As $W$ is a canonical $\phi$-coordinated
$V_{W}$-module, $W$ is a $\phi$-coordinated
$V_{q}(\beta\gamma)$-module.

On the other hand, assume that $(W,Y_{W})$ is a $\phi$-coordinated
$V_{q}(\beta\gamma)$-module. In view of Propositions
\ref{palgebra-module} and \ref{pwqva-phimodule}, we have
\begin{eqnarray*}
&&Y_{W}(\hat{\beta},x_{1})Y_{W}(\hat{\beta},x_{2})
=\iota_{x_{2},x_{1}}\left(\frac{x_{1}-qx_{2}}{qx_{1}-x_{2}}\right)
Y_{W}(\hat{\beta},x_{2})Y_{W}(\hat{\beta},x_{1}),\\
&&Y_{W}(\hat{\gamma},x_{1})Y_{W}(\hat{\gamma},x_{2})
=\iota_{x_{2},x_{1}}\left(\frac{x_{1}-qx_{2}}{qx_{1}-x_{2}}\right)
Y_{W}(\hat{\gamma},x_{2})Y_{W}(\hat{\gamma},x_{1}),\\
&&Y_{W}(\hat{\beta},x_{1})Y_{W}(\hat{\beta},x_{2})
-\iota_{x_{2},x_{1}}\left(\frac{qx_{1}-x_{2}}{x_{1}-qx_{2}}\right)
Y_{W}(\hat{\gamma},x_{2})Y_{W}(\hat{\beta},x_{1})\\
&&\hspace{1cm}=\Res_{x_{0}}x_{1}^{-1}\delta\left(\frac{x_{2}e^{x_{0}}}{x_{1}}\right)
x_{2}e^{x_{0}}Y_{W}(Y(\hat{\beta},x_{0})\hat{\gamma},x_{2})\\
 &&\hspace{1cm}=\delta\left(\frac{x_{2}}{x_{1}}\right),
\end{eqnarray*}
where we are using the relation
$$\hat{\beta}_{n}\hat{\gamma}=\delta_{n,0}{\bf 1}\ \ \ \mbox{ for
}n\ge 0.$$ Thus, $W$ is a restricted
$A_{q}^{trig}(\beta\gamma)$-module with $\tilde{\beta}(z)$ and
$\tilde{\gamma}(x)$ acting as $Y_{W}(\hat{\beta},x)$ and
$Y_{W}(\hat{\gamma},x)$, respectively.
\end{proof}


\begin{thebibliography}{FLM}

\bibitem [BK]{bk}
B. Bakalov and V. Kac, Field algebras, {\em Internat. Math. Res.
Notices} {\bf 3} (2003) 123-159.

\bibitem[Bo]{boch}
S. Bochner, Formal Lie groups, {\em Ann. of Math.} {\bf 47} (1946)
192-201.

\bibitem[Bor]{b-Gva}
R. E. Borcherds, Vertex algebras, in ``{\em Topological Field
Theory, Primitive Forms and Related Topics}'' (Kyoto, 1996), edited
by M. Kashiwara, A. Matsuo, K. Saito and I. Satake, Progress in
Math., Vol. 160, Birkh\"auser, Boston, 1998, 35-77.

\bibitem[Dr]{drinfeld}
V. G. Drinfeld, A new realization of Yangians and quantized affine
algebras, {\em Soviet Math. Dokl.} {\bf 36} (1988) 212-216.

\bibitem[EFK]{efk}
P. Etingof, I. Frenkel and A. Kirillov Jr., {\em Lectures on
Representation Theory and Knizhnik-Zamolodchikov Equations},
Mathematical Surveys and Monographs, Vol. 58, Amer. Math. Soc.,
Providence, 1998.

\bibitem[EK]{ek}
P. Etingof and D. Kazhdan, Quantization of Lie bialgebras, V, {\em
Selecta Mathematica (New Series)} {\bf 6} (2000) 105-130.

\bibitem[Fad]{fad}
L. Faddeev, Quantum completely integrable models in field theory,
{\em Soviet Sci. Rev., Ser. C: Math. Phys. Rev.} {\bf 1}, Hawood
Academic Publ., 1990, pp. 107-155.

\bibitem [FHL]{fhl}
I. B. Frenkel, Y.-Z. Huang, and J. Lepowsky, On axiomatic approaches
to vertex operator algebras and modules, Memoirs Amer. Math. Soc.
{\bf 104}, 1993.

\bibitem [FJ]{fj}
I. B. Frenkel and N.-H. Jing, Vertex representations of quantum
affine algebras, {\em Proc. Natl. Acad. Sci. USA} {\bf 85} (1988)
9373-9377.

\bibitem[FLM]{flm}
I. Frenkel, J. Lepowsky and A. Meurman, {\it Vertex Operator
Algebras and the Monster}, Pure and Appl. Math., {\bf Vol. 134},
Academic Press, Boston, 1988.

\bibitem[FZ]{fz}
I. Frenkel and Y.-C. Zhu, Vertex operator algebras associated to
representations of affine and Virasoro algebras, {\it Duke Math. J.}
{\bf 66} (1992) 123-168.

\bibitem [KL]{kl}
M. Karel and H.-S. Li, Some quantum vertex algebras of
Zamolodchikov-Faddeev type, {\em Commun. Contemp. Math.} {\bf 11}
(2009) 829-863.

\bibitem[LL]{ll}
J. Lepowsky and H.-S. Li, {\em Introduction to Vertex Operator
Algebras and Their Representations}, Progress in Math.  {\bf 227},
Birkh\"auser, Boston, 2003.

\bibitem[Li1]{li-local}
H.-S. Li, Local systems of vertex operators, vertex superalgebras
and modules, {\em J. Pure Appl. Alg.} {\bf 109} (1996) 143-195;
hep-th/9406185.

\bibitem[Li2]{li-g1}
H.-S. Li, Axiomatic $G_{1}$-vertex algebras, {\em Commun. Contemp.
Math.} {\bf 5} (2003) 281-327.

\bibitem[Li3]{li-qva1}
H.-S. Li, Nonlocal vertex algebras generated by formal vertex
operators, {\em Selecta Mathematica (New Series)} {\bf 11} (2005)
349-397.

\bibitem[Li4]{li-qva2}
H.-S. Li, Constructing quantum vertex algebras, {\em International
Journal of Mathematics} {\bf 17} (2006) 441-476.

\bibitem[Li5]{li-gamma}
H.-S. Li, A new construction of vertex algebras and quasi modules
for vertex algebras, {\em Advances in Math.} {\bf 202} (2006)
232-286.

\bibitem[Li6]{li-infinity}
H.-S. Li, Modular-at-infinity for quantum vertex algebras, {\em
Commun. Math. Phys.} {\bf 282} (2008) 819-864.

\bibitem[Li7]{li-hqva}
H.-S. Li, $\hbar$-adic quantum vertex algebras and their modules,
{\em Commun. Math. Phys.} {\bf 296} (2010) 475-523.

\bibitem[Li8]{li-tqva}
H.-S. Li, Quantum vertex ${\Bbb F}((t))$-algebras and their modules,
preprint, 2009; arXiv:0903.0186 [math.QA].

\bibitem[LTW]{ltw}
H.-S. Li, Shaobin Tan and Qing Wang, Twisted modules for quantum
vertex algebras, {\em J. Pure Applied Algebra} {\bf 214} (2010)
201-220.

\bibitem[ZZ]{zam}
A. B. Zamolodchikov and Al. B. Zamolodchikov, Factorized
$S$-matrices in two dimensionals as the exact solutions of certain
relativistic quantum field theory models, {\em Ann. of Physics} {\bf
120} (1979) 253-291.

\end{thebibliography}
\end{document}